\documentclass[12pt]{article}
\usepackage[utf8]{inputenc}

\usepackage[style = apa]{biblatex}
\DeclareLanguageMapping{english}{english-apa}
\let\cite\parencite

\let\vec\mathbf
\usepackage{lipsum}
\usepackage{amsfonts}
\usepackage{geometry}
\usepackage{graphicx}
\usepackage{epstopdf}
\usepackage{xcolor}
\usepackage{mathtools}
\usepackage{amssymb}
\usepackage{esint}
\usepackage{bm}
\usepackage{array}
\usepackage{float}
\usepackage{caption}
\usepackage{subcaption}
\usepackage{enumerate}
\usepackage{algpseudocode}
\usepackage{algorithm}
\usepackage[hidelinks]{hyperref}
\usepackage{cleveref}
\crefname{equation}{}{}
\crefname{figure}{Figure}{}

\let\vec\mathbf

\algblock{ParFor}{EndParFor}
\algnewcommand\algorithmicparfor{\textbf{parfor}}
\algnewcommand\algorithmicpardo{\textbf{do}}
\algnewcommand\algorithmicendparfor{\textbf{end\ parfor}}
\algrenewtext{ParFor}[1]{\algorithmicparfor\ #1\ \algorithmicpardo}
\algrenewtext{EndParFor}{\algorithmicendparfor}

\newcommand{\dt}[0]{\hspace{0.75mm}\text{d} t }

\newcommand{\dbydt}[1]{\frac{\text{d} {#1}}{\dt}}

\title{A mean correction for improved phase-averaging accuracy in oscillatory, multiscale, differential equations}

\author{Timothy C. Andrews \thanks{University of Exeter (ta440@exeter.ac.uk, b.wingate@exeter.ac.uk) }
\and Beth A. Wingate\footnotemark[1]}

\date{November 5, 2024}

\begin{document}

\maketitle

\begin{abstract}
This paper introduces a new algorithm to improve the accuracy of numerical phase-averaging in oscillatory, multiscale, differential equations. Phase-averaging is a timestepping method which averages a mapped variable to remove highly oscillatory linear terms from the differential equation. This retains the main contribution of fast waves on the low frequencies without explicitly resolving the rapid oscillations. However, this comes at the cost of introducing an averaging error. To offset this, we propose a modified mapping that includes a mean correction term encoding an average measure of the nonlinear interactions. This mapping was introduced in [\textit{Commun Nonlinear Sci Numer Simulat 71 (2019) 1–21}] for weak nonlinearity and relied on classical time-averaging, which leaves only the zero frequencies. Our algorithm instead considers mean corrected phase-averaging when 1) the nonlinearity is not weak but the linear oscillations are fast and 2) finite averaging windows are applied via a smooth kernel, which has the advantage of retaining low frequencies whilst still eliminating the fastest oscillations. In particular, we introduce a local mean correction that combines the concepts of a mean correction and finite averaging; this retains low-frequency components in the mean correction that are removed with classical time-averaging. We show that the new timestepping algorithm reduces phase errors in the mapped variable for the swinging spring ODE in various dynamical configurations. We also show accuracy improvements with a local mean correction compared to standard phase-averaging in the one-dimensional rotating shallow water equations, a useful test case for weather and climate applications. 
\end{abstract}

\section{Introduction}
\label{sec:intro}
The ability to take larger timesteps, whilst retaining sufficient accuracy, is of great importance to reduce computational expense and time-to-solution. A particular motivation is weather and climate modelling, where the governing equations are highly oscillatory. In previous decades, research into semi-Lagrangian methods \cite{temperton1987_SLSI,staniforth_SL_review,SL_highres_ECMWF} and implicit-explicit (IMEX) timesteppers \cite{giraldo2013_hevi, weller2013_RK_IMEX, bao2015_hevi} enabled larger timesteps for global weather forecasting models without sacrificing accuracy \cite{Mengaldo_time_integration_NWP}. An alternative approach that we will examine is the phase-averaging method. This has the potential for further reductions in time-to-solution through time parallelism \cite{Haut_Wingate_2014,hiroe_colin}. This paper proposes a new algorithm to improve the accuracy of this phase-averaging method. \par
Phase-averaging is a numerical technique that uses mapping and averaging to smooth out fast linear oscillations, whilst retaining the lower frequencies. First, a mapping to a modulation variable enables a time evolution without the leading-order oscillations, which can increase the allowable timestep size of explicit methods. Next, a phase-average of the modulation variable timestepping equation is applied to smooth the higher-order oscillations. This further reduces the numerical stiffness in the equations to enable larger timesteps, but comes at the cost of introducing an averaging error. We seek to increase the accuracy of the phase-averaged method presented in \cite{peddle_haut_wingate_2019} by including a mean correction term, which encodes an average of the nonlinearity, into the modulation variable mapping. This technique was introduced by Tao in \cite{Tao_simp_imp} for weakly nonlinear problems using classical averaging techniques \cite{sanders2007averaging}. Our algorithm in this paper can be considered a generalisation of Tao's method for cases where the equation is not weakly nonlinear, but oscillating and multiscale, and finite averaging windows are applied via a smooth kernel. 
\par
Classical averaging applied to oscillatory systems, as in \cite{Tao_simp_imp}, retains only the zero frequency (non-oscillatory) components in the mapped variable. A finite averaging window instead permits small, but nonzero, frequencies, which are important for dynamics over long time periods, e.g. \cite{Smith_Lee,chen2005resonant}. The benefits of a finite averaging window are discussed in detail in \cite{Adam}; here, it was shown that the error bound of this approach is a function of the timestep size, $\Delta T$, the averaging window, $\eta$, and a parameter, $\epsilon$, which represents the speed of the fast linear oscillations. Advantageously, the error bound for the finite phase-average not only holds in the limit of $\epsilon$ small, as in classical averaging, but also for finite values of $\epsilon$. 
\par
We consider the application of phase-averaging to multiscale initial value problems of the form 
\begin{equation}
\frac{\partial \vec{U}}{\partial t} + \frac{1}{\epsilon} \mathcal{L} \vec{U} = \mathcal{N}(\vec{U}), ~\vec{U}(0)=\vec{U}_0, ~\epsilon \in (0,1], ~\epsilon \in \mathbb{R},
\label{eq:standard_eq}
\end{equation}

\noindent where $\vec{U}(t)$ is a (spatial) vector-valued function, $\vec{U}(t) = (u_1(t,\cdot),u_2(t,\cdot),...)$, $\mathcal{L}$ has purely imaginary (oscillatory) eigenvalues, and $\mathcal{N}(\vec{U})$ is an autonomous nonlinear term with arbitrary order polynomial components. In this paper, we term (\ref{eq:standard_eq}) the standard equation, which is an ordinary differential equation (ODE) or a partial differential equation (PDE) that has been discretised in space to form a system of ODEs. The parameter $\epsilon$, which is assumed to be small, measures how rapid the linear oscillations are.  The need to sufficiently resolve the fastest of these frequencies on an $\mathcal{O}(\epsilon^{-1})$ timescale constrains the maximum allowable timestep of explicit timestepping methods in (\ref{eq:standard_eq}). This restriction to small timesteps is also true for implicit methods if accuracy is required.  \par
The phase-averaging approach first introduces a matrix exponential mapping to a \textit{modulation variable}, $\vec{V}$:
\begin{equation}
\vec{U}(t) = e^{-\frac{t}{\epsilon} \mathcal{L} } \vec{V}(t).
\label{eq:modvar}
\end{equation}

\noindent Using the modulation variable in the standard equation (\ref{eq:standard_eq}) leads to a new system to timestep,
\begin{equation}
\frac{\partial \vec{V}}{\partial t} = e^{\frac{t}{\epsilon} \mathcal{L}} \mathcal{N} \left( e^{- \frac{t}{\epsilon} \mathcal{L}} \vec{V} \right), \ \ \  \vec{V}(0) = \vec{U}_0.
\label{eq:modvar_system}
\end{equation} 

\noindent The standard equation can be recovered by applying the inverse of the matrix exponential to (\ref{eq:modvar_system}). Advantageously, the modulation variable system (\ref{eq:modvar_system}) has no linear terms oscillating on an $\mathcal{O}(\epsilon^{-1})$ timescale, so there is a reduction in numerical stiffness compared to (\ref{eq:standard_eq}). The use of mapping to a modulation variable equation is standard practice in the method of averaging \cite{sanders2007averaging,BoMi1961} and has been used to solve stiff PDEs in different fields, including optics \cite{interaction_picture_GNLSE} and geophysical fluid flows  \cite{Jones_splitting_SWEs,Embid_1998,BaMaNi1998,BaMaNiZh1997, KlMa1981}. \par
Taking a time derivative of (\ref{eq:modvar_system}) shows that there are still oscillations in the modulation variable space that can impact large timestep stability and accuracy. Therefore, phase-averaging is used to smooth these small amplitude oscillations and further reduce the numerical stiffness. We facilitate this by introducing a phase variable, $s$, into the linear mapping of (\ref{eq:modvar}) to specify a time-translation, as in \cite{bauer_colin_spring,Adam},
\begin{equation}
    \vec{U}(t,s) = e^{- \frac{t + s}{\epsilon} \mathcal{L} } \vec{V}(t,s),
    \label{eq:modvar_with_s}
\end{equation}

\noindent so that the phase variable appears in the timestepping equation (\ref{eq:modvar_system}),
\begin{equation}
\frac{\partial \vec{V} (t,s) }{\partial t} = e^{\frac{t + s}{\epsilon} \mathcal{L} } \mathcal{N} \left( e^{- \frac{t + s}{\epsilon} \mathcal{L} } \vec{V} \right).
\label{eq:modvar_system_with_s}
\end{equation} 

\noindent We then integrate over this phase variable to compute an average, denoted as $\left< \cdot \right>_s$, leaving only time-dependence in the phase-averaged modulation variable,
\begin{equation}
\frac{\partial \vec{\overline{V}}(t)}{\partial t} = \left< e^{\frac{t + s}{\epsilon} \mathcal{L} } \mathcal{N} \left( e^{- \frac{t + s}{\epsilon} \mathcal{L} } \vec{\overline{V}} \right) \right>_s, ~\vec{\overline{V}}(0) = \vec{U}_0.
\label{eq:modvar_angle_bracket_ave}
\end{equation}

The method of averaging has a long history and a review can be found in \cite{sanders2007averaging}. It has been used to obtain numerical solutions of oscillatory differential equations in many different forms (see \cite{Tao_simp_imp} for an overview). It was first proposed to solve systems like (\ref{eq:standard_eq}) using a predictor-corrector type parallel-in-time method by \cite{Haut_Wingate_2014}. Their phase-averaging was classical, which corresponds to $s \in (-\infty,\infty)$ and the removal of all non-zero frequencies in the timestepping equation. In \cite{peddle_haut_wingate_2019} this was extended to the more physically relevant case of a finite oscillation speed, $\epsilon$. This phase-average is computed over only a few oscillations using a finite window length, $\eta$, with $s \in ( -\eta / 2, \eta / 2 )$.  The algorithm for finite phase-averaging has been extended to a finite element implementation of the rotating shallow water equations on the sphere \cite{hiroe_colin} and a multi-level parareal method that uses different averaging windows on each level \cite{rosemeier2024multilevel}. \par
The averaging error introduced by $\left< \cdot \right>_s$ phase-average operation is balanced with the reduced numerical stiffness which allows a larger timestep in the averaged system. This leads to the possibility of choosing an $(\eta,\Delta t)$ pair that achieves a minimum error for a given $\epsilon$ \cite{Adam}. Our new algorithm in this paper seeks to reduce this minimum error by introducing an additional term into the modulation variable mapping (\ref{eq:modvar}). This new term quantifies a local mean of the oscillations in the modulation variable space. Its use within the mapping ensures that less information is lost when phase-averaging the modulation variable equation of (\ref{eq:modvar_system_with_s}). The extension of phase-averaging with a mean correction can be considered an alternative to \cite{bauer_colin_spring}, where higher-order versions of the phase-averaging are instead used to improve the accuracy. 
\par
For the new method we will investigate in this paper, we define a modified mapping to a different modulation variable, $\vec{W}$, 
\begin{equation}
\vec{U}(t, s) = e^{-\frac{t + s}{\epsilon} \mathcal{L}} \vec{W}(t,s) + \epsilon \mathcal{L}^{-1} \vec{C},
\label{eq:meancor_modvar}
\end{equation}

\noindent where the additional $\epsilon \mathcal{L}^{-1} \vec{C}$ term enables the mapped oscillations to be about some non-zero mean value. The inclusion of the phase variable, $s$, is an extension of Tao's mapping in \cite{Tao_simp_imp} that allows the classical time-average to be replaced with a finite phase-average, resulting in:
\begin{equation}
\begin{split}
\frac{\partial \vec{\overline{W}}(t)}{\partial t} = \left< e^{ \frac{t+s}{\epsilon} \mathcal{L}} \left[ \mathcal{N} \left( e^{- \frac{t+s}{\epsilon} \mathcal{L}} \vec{\overline{W}}(t) + \epsilon \mathcal{L}^{-1} \vec{C} \right) - \mathcal{L} \mathcal{L}^{-1} \vec{C} \right] \right>_s, \\
\vec{\overline{W}}(0) = \vec{U_0} - \epsilon \mathcal{L}^{-1} \vec{C}(\vec{\overline{W}}(0)).
\label{eq:meancor_modvar_angle_bracket_ave}
\end{split}
\end{equation}

A key idea in this paper is the creation of a local mean correction that is computed using a finite averaging window. Analogously to a classical time-average of the timestepping equation (\ref{eq:modvar_angle_bracket_ave}), the classical mean correction of \cite{Tao_simp_imp} only includes interactions in the nonlinearity that have a frequency of zero. By contrast, the local mean correction retains low-frequency components to enable greater accuracy with a finite phase-average of the timestepping equation. In two numerical PDE experiments (sections \ref{sec:kg} and \ref{sec:1drswe}), we will identify best finite windows for computing the local mean correction; this choice is dependent on the speed of the linear oscillations.
\par
In this paper, we will explore potential accuracy improvements of using (\ref{eq:meancor_modvar_angle_bracket_ave}) over (\ref{eq:modvar_angle_bracket_ave}) for timestepping in oscillatory, multiscale, PDEs. Section \ref{sec:phaseaverage} will delve into the phase-averaged and mean corrected methods, including pseudocode for the new algorithm. We will highlight two key differences between our method for multiscale equations and the previous application in \cite{Tao_simp_imp}: i) a finite phase-average is used instead of a classical time-average (section \ref{subsec:phase_average}), and ii) a \textit{local} mean correction using a finite time interval ($\vec{C}$) can be used instead of the classically time-averaged version ($\vec{C}_{\infty}$) (section \ref{subsec:meancor_method}). Intuition around accuracy improvements with the mean correction will be provided in section \ref{subsec:simple_meancor_example}, where we consider a simple ODE with an analytical solution. We will then explore the potential accuracy improvements of the new algorithm in three numerical experiments (sections \ref{sec:ss}-\ref{sec:1drswe}). Section \ref{sec:ss} will present the swinging spring ODE \cite{swing_spring}, where we will compute a classical mean correction and see it correcting phase errors in the modulation variable space. The first PDE example in section \ref{sec:kg} is a Klein-Gordon-type (KG-type) PDE, which is both weakly nonlinear and multiscale. Section \ref{sec:1drswe} will explore accuracy improvements in the one-dimensional rotating shallow water equations (RSWEs), a prototype PDE for those used in weather and climate applications. In section \ref{sec:discussion} we will discuss future improvements that can be made to the locally mean corrected method and list further applications that could be investigated.

\section{Phase-averaging and the mean correction}
\label{sec:phaseaverage}

\subsection{Discussion of phase-averaging}
\label{subsec:phase_average}
Finite phase-averaging differs from a classical time-average in that it is computed over a select number of oscillations. The phase-averaging technique introduces a dummy phase variable, $s$, into the linear mapping and integrates over this variable with respect to a weighted kernel function. The range of phase shifts used in the average are those contained within the averaging window of length $\eta$. We define the $\left< \cdot \right>_s $ average as
\begin{equation}
\langle \phi \rangle_s = \int_{s = - \eta / 2}^{\eta / 2} \mathcal{K}(s) ~\phi  ~\text{d}s,
\end{equation}

\noindent where $\mathcal{K}(s)$ denotes the averaging kernel function. Substituting this into {(\ref{eq:meancor_modvar_angle_bracket_ave})} obtains:
\begin{equation}
\frac{\partial \vec{\overline{V}} (t)}{\partial t} = \int_{s = - \eta / 2}^{\eta / 2} \mathcal{K}(s) e^{ \frac{t+s}{\epsilon} \mathcal{L}} \mathcal{N} \left( e^{- \frac{t+s}{\epsilon} \mathcal{L}} \vec{\overline{V}}(t) \right) \text{d} s.
\label{eq:phase_ave_modvar}
\end{equation}

\noindent This approach of numerical averaging with a kernel is analogous to forming the effective system in heterogeneous multiscale methods \cite{engquist2005_HMM}. \par
There are three regimes of averaging window length that could be used in (\ref{eq:phase_ave_modvar}) --- case iii) is of specific interest in this paper:
 \vspace{.1in}
\begin{enumerate}[i)]
    \item The limit of $\eta = 0$ implies no averaging, returning to the usual modulation variable system (\ref{eq:modvar_system}). This equation remains constrained by small amplitude oscillations.
    \item The other limit as $\eta \rightarrow \infty$ simulates the limit of $\epsilon \rightarrow 0$. This solution has a complete cancellation of oscillations; any parts of the solution with periodicity, and thus phase dependence, are driven to zero. As such, only non-oscillatory components contribute to the leading-order solution \cite{Schochet}. \cite{Haut_Wingate_2014} showed that the dynamics in this limit can be replicated with a sufficiently large $\eta$.
    \item An averaging window with minimum error, $\eta^*$, at a fixed timestep of $\Delta t$, was shown to exist in finite ($\epsilon$) regimes \cite{peddle_haut_wingate_2019}. Insufficient averaging $(\eta < \eta^*)$ means that the timestepping equations remain too stiff, which leaves a significant timestepping error. Excessive averaging $(\eta > \eta^*)$ removes dynamical information contained in the oscillations, so the averaging error dominates. The combination of timestepping error (decreasing with $\eta$) and averaging error (increasing with $\eta$) is smallest at the averaging window of $\eta^*$. Visualisations of phase-averaged solution errors against $\eta$ allow the identification of $\eta^*$; these will be termed \textit{Peddle plots} (see Figure \ref{fig:swing_spring_peddle} for an example in the swinging spring ODE, or \cite{Adam} for others). As $\eta^*$ depends on the timestep size, errors are computed over a range of normalised averaging windows, $\zeta = \eta / \Delta t$.
\end{enumerate}
\vspace{.1in}
\par
An exponential bump function is selected for the averaging kernel to ensure sufficient smoothness when approximating the $\epsilon \rightarrow 0$ limit with a compact support \cite{Haut_Wingate_2014}. The kernel has a maximum value at no phase shift ($s=0$) and decays to zero at the boundaries of the averaging interval,
\begin{equation}
\mathcal{K}(s)
=
\begin{cases}
A \exp \left( \frac{-\eta}{\Gamma (0.5 \eta + s) (0.5 \eta - s)} \right), s \in \left( -\frac{\eta}{2},\frac{\eta}{2} \right), \\
0, |s| \geq \frac{\eta}{2},
\end{cases}
\label{eq:kernel_func}
\end{equation}

\noindent where $A$ is a normalisation constant to ensure that $\int_{s=-\infty}^{\infty} \mathcal{K}(s) \  \text{d} s = 1$. A kernel decay rate, $\Gamma$, controls the exponential's shape over the compact support. Values of $\Gamma=1$ \cite{peddle_haut_wingate_2019,hiroe_colin}  and $\Gamma=4$ \cite{Adam} have been used in prior applications --- we use the latter for the results in this paper. 
\par
The numerical implementation can be expressed as a sum of evaluations at $K_s$ equispaced phase shifts, 
\begin{equation}
    \frac{\partial \vec{\overline{V}} (t)}{\partial t} = \sum_{k=1}^{K_s} \mathcal{K}(s_k) e^{ \frac{t+s_k}{\epsilon} \mathcal{L}} \mathcal{N} \left( e^{- \frac{t+s_k}{\epsilon} \mathcal{L}} \vec{\overline{V}}(t) \right), ~s_k = \left( \frac{2k-1}{2K} - \frac{1}{2} \right) \eta,
    \label{eq:phase_ave_discrete}
\end{equation}

\noindent where the discrete kernel values are normalised to ensure that $\sum_{k=1}^{K_s} \mathcal{K}(s_k) = 1$. 

\par
Convergence of the numerical phase-average requires a sufficient number of discrete kernel points, $K_s$. For PDE systems, this scales with the fastest resolvable frequency and the averaging window length. The following rule is used to compute $K_s$:
\begin{equation}
    K_s = \Big\lceil P \eta \frac{\omega_{\text{max}}}{2 \pi \epsilon} \Big\rceil,
    \label{eq:K_rule}
\end{equation}

\noindent where $\lceil . \rceil$ denotes the ceiling function, $\omega_{\text{max}}$ is the fastest resolvable frequency of $\mathcal{L}$, and $P$ is the number of samples within each period of this fastest frequency. To have sufficient sampling of the fastest oscillations, a value of at least $P = 3.5$ is recommended for pseudospectral methods \cite{boyd2001chebyshev}. Following \cite{hiroe_colin}, we will use $P = 4$. As a proportionally larger number of points per oscillation may be required when sampling fewer oscillations \cite{boyd2001chebyshev}, a minimum $K_s$ value is also specified. \par
A key feature of the phase-average in (\ref{eq:phase_ave_modvar}) is that the oscillations in the matrix exponential are averaged at a \textit{fixed} value of the solution, $\vec{\overline{V}}(t)$. This only requires phase-shifting the linear operator of $\mathcal{L}$. As the computations at each $s$ value have no feedback on each other, the phase-average can be computed in parallel using $K_s$ processors \cite{Haut_Wingate_2014}.

\subsection{Mean corrected phase-averaging}
\label{subsec:meancor_method}
With standard phase-averaging, the mapping to the modulation variable (\ref{eq:modvar}) effectively assumes that the oscillations are around a state of zero. The idea of \cite{Tao_simp_imp} is that we can improve this transform by allowing for oscillations about some non-zero mean value. To derive this additional term, we assume that the nonlinearity is slowly varying and can be approximated as a constant \textit{over one timestep}: $\mathcal{N}(\vec{W}) \sim \vec{C} (\vec{W} (t))$, where $\vec{W}$ is a modulation variable and $\vec{C}$ is a constant that depends on the current solution value of $\vec{W}(t)$. We then seek a mean component of the modulation variable, $\vec{W}_M$, by substituting this into our standard equation (\ref{eq:standard_eq}) with the constant nonlinearity:
\begin{align}
    \frac{\partial \vec{W}_M}{\partial t} + \epsilon^{-1} \mathcal{L} \vec{W}_M &= \vec{C} \nonumber \\
    \rightarrow \vec{W}_M &= \epsilon \mathcal{L}^{-1} \vec{C}.
\end{align}

We then add this mean term to the original modulation variable definition to obtain the mean correction mapping (\ref{eq:meancor_modvar}), which we restate here without the phase shift:
\begin{equation}
\vec{U}(t) = e^{-\frac{t}{\epsilon} \mathcal{L}} \vec{W}(t) + \epsilon \mathcal{L}^{-1} \vec{C}
\end{equation}

\noindent $\vec{C}$ quantifies an average effect of nonlinear oscillations of the modulation variable,
\begin{equation}
\vec{C} =  \left< \mathcal{N} \left( e^{-\frac{t + r}{\epsilon} \mathcal{L}} \vec{W}(t) \right) \right>_r,
\label{eq:C_general}
\end{equation}
where $r$ is an analogous phase variable to $s$ and $\left< \cdot \right>_r$ denotes the averaging operation to compute the mean correction. 
\par
The application for weakly nonlinear systems in \cite{Tao_simp_imp} used a classical time-average to compute the mean correction:
\begin{equation}
\vec{C}_{\infty} (\vec{W}(t)) =  \lim_{T_C \rightarrow \infty} \frac{1}{T_C} \int_{r = 0 }^{T_C} \mathcal{N} \left( e^{-\frac{t + r}{\epsilon} \mathcal{L}} \vec{W}(t) \right) \text{d} r .
\label{eq:C_inf}
\end{equation}

\noindent In some situations, it is feasible to compute a closed-form expression for $\vec{C}_{\infty}$ using (\ref{eq:C_inf}); we call this the classical mean correction. A classical mean correction will be used in the ODE system of the swinging spring (section 3) and the KG-type system (section 4); this is made tractable by our use of a pseudospectral spatial discretisation in this work. However, our primary focus is on the new idea of using a finite averaging window in (\ref{eq:C_general}) to compute $\vec{C}$. The corresponding \textit{local} mean correction has two key advantages over the classical mean correction for use with finite phase-averaging:
\begin{itemize}
    \item It retains low-frequency components of nonlinear oscillation that are lost when taking $T_C \rightarrow \infty$.
    \item It can computed algorithmically for any numerical model as opposed to being derived for a specific equation and spatial discretisation choice.
\end{itemize}

The local mean correction is computed using an analogous phase-averaging computation for $\left< . \right>_r$ as is used for $\left< . \right>_s$,
\begin{equation}
\vec{C} (\vec{W}(t), t) =  \int_{r = -\eta_C / 2 }^{\eta_C/2} \mathcal{K}_C(r) \ \mathcal{N} \left( e^{-\frac{t + r}{\epsilon} \mathcal{L}} \vec{W}(t) \right) \ \text{d} r ,
\label{eq:C_finite}
\end{equation}

\noindent where $\mathcal{K}_C$ is the mean correction averaging kernel of length $\eta_C$. We select the same functional form for $\mathcal{K}_C(r)$ as is defined for $\mathcal{K}(s)$ in (\ref{eq:kernel_func}). The required $K_r$ number of points in the mean correction kernel is computed using (\ref{eq:K_rule}) with $\eta_C$ in place of $\eta$. Algorithm \ref{alg:num_C} provides pseudocode for the computation of the local mean correction. 

\begin{algorithm}
\caption{The Local Mean Correction}
\label{alg:num_C}
C$(\vec{W}^n, t^n):$
\begin{algorithmic}[1]
    \ParFor{$k \in \{ 1, 2, ..., K_s \}$ }
        \State $\vec{U}_k \gets  \exp (- ((t^n + r_k) / \epsilon) \mathcal{L}) \vec{W}^n$
        \State $\mathcal{N}_k \gets \mathcal{N}(\vec{U}_k)$
            \EndParFor
    \State $\vec{C} \gets \sum_{k=1}^{K_s} \mathcal{K}_C(r_k) \mathcal{N}_k$ 
\State Return $\vec{C}$
\end{algorithmic}
\end{algorithm}

As with the $\left< \cdot \right>_s$ phase-averaging procedure, the local mean correction is computed at a fixed solution state, $\vec{W}(t)$, so is parallelisable. The locally mean corrected method requires selecting a suitable window pairing of $(\eta, \eta_C)$. Analogously to an $\eta^*$ minimising the error of the $\left< . \right>_s$ phase-average, the existence of a finite $\eta_C^*$ which minimises the error of the locally mean corrected method is observed in sections \ref{sec:kg} and \ref{sec:1drswe}.
\par
We now apply the new modulation variable transformation (\ref{eq:meancor_modvar}), using either the classical or local mean correction, to the standard form equation (\ref{eq:standard_eq}),

\begin{align}
    \frac{\partial}{\partial t} \left( e^{-\frac{t + s}{\epsilon} \mathcal{L}} \vec{W} + \epsilon \mathcal{L}^{-1} \vec{C} \right) + \frac{1}{\epsilon} \mathcal{L} \left( e^{-\frac{t + s}{\epsilon} \mathcal{L}} \vec{W} + \epsilon \mathcal{L}^{-1} \vec{C} \right) = \mathcal{N}, \nonumber \\ 
    \rightarrow -\frac{1}{\epsilon} \mathcal{L} e^{-\frac{t + s}{\epsilon} \mathcal{L}} \vec{W} + e^{-\frac{t + s}{\epsilon} \mathcal{L}}  \frac{\partial \vec{W}}{\partial t} + \epsilon \frac{\partial }{\partial t} (\mathcal{L}^{-1} \vec{C}) + \frac{1}{\epsilon} \mathcal{L} e^{-\frac{t + s}{\epsilon} \mathcal{L}} \vec{W}   + \mathcal{L} \mathcal{L}^{-1} \vec{C} = \mathcal{N}, \nonumber \\
    \rightarrow e^{-\frac{t + s}{\epsilon} \mathcal{L}} \frac{\partial \vec{W}}{\partial t} + \epsilon \frac{\partial }{\partial t} (\mathcal{L}^{-1} \vec{C}) = \mathcal{N} - \mathcal{L} \mathcal{L}^{-1} \vec{C}. \nonumber
\end{align}

\noindent Applying the chain rule to the $\epsilon \frac{\partial }{\partial t} (\mathcal{L}^{-1} \vec{C})$ term,

\begin{equation}
    \left( e^{-\frac{t + s}{\epsilon} \mathcal{L}} + \epsilon \frac{\partial}{\partial \vec{W}} (\mathcal{L}^{-1} \vec{C}) \right)   \frac{\partial \vec{W}}{\partial t} = \mathcal{N} - \mathcal{L} \mathcal{L}^{-1} \vec{C}.
    \label{eq:meancor_modvar_with_eps_d_dw}
\end{equation}

\noindent Following \cite{Tao_simp_imp} we assume that we are in a regime of $\epsilon << 1$, which in our multiscale equations corresponds to sufficiently fast linear oscillations. As such, the $\epsilon \frac{\partial}{\partial \vec{W}} (\mathcal{L}^{-1} \vec{C})$ in (\ref{eq:meancor_modvar_with_eps_d_dw}) term can be neglected at $\mathcal{O} (1)$. The validity of this assumption will be revisited in sections \ref{sec:kg} and \ref{sec:1drswe} in numerical PDE experiments at moderate values of $\epsilon$. \par
Using the assumption of fast oscillations leaves the following time evolution equation for the mean corrected modulation variable,
\begin{equation}
\frac{\partial \vec{W} (t,s)}{\partial t} = e^{\frac{t + s}{\epsilon} \mathcal{L}} \left[ \mathcal{N} \left( e^{-\frac{t + s}{\epsilon} \mathcal{L}} \vec{W} + \epsilon \mathcal{L}^{-1} \vec{C} \right) - \mathcal{L} \mathcal{L}^{-1} \vec{C}  \right].
\end{equation}
\par

\noindent This equation is phase-averaged over a finite time window to give the following numerical system for $\overline{\vec{W}}$,
\begin{equation}
\begin{split}
\frac{\partial \vec{\overline{W}} (t)}{\partial t} = \int_{s = -\eta / 2}^{\eta / 2} \mathcal{K}(s) e^{ \frac{t+s}{\epsilon} \mathcal{L}} \left[ \mathcal{N} \left( e^{ - \frac{t+s}{\epsilon} \mathcal{L}} \vec{\overline{W}} + \epsilon \mathcal{L}^{-1} \vec{C}  \right) - \mathcal{L} \mathcal{L}^{-1} \vec{C} \right]  \text{d}s,   \\
\overline{\vec{W}}(0) = \vec{U_0} - \epsilon \mathcal{L}^{-1} \vec{C} (\overline{\vec{W}}(0)).
\end{split}
\label{eq:aved_meancor_modvar_system}
\end{equation}

The new algorithm that solves (\ref{eq:aved_meancor_modvar_system}) is given in Algorithm \ref{alg:meancor_alg}. This is presented in a format where the averaging procedures of $\left< \cdot \right>_s$ and $\left< \cdot \right>_r$ can be computed in parallel (using parfor). As this paper is concerned with exploring potential accuracy improvements, not parallelisation, our results will use a serial implementation. 
\par
Many systems have non-singular linear matrices, so the $- \mathcal{L} \mathcal{L}^{-1} \vec{C} $ term in (\ref{eq:aved_meancor_modvar_system}) often simplifies to $- \vec{C} $. However, this reduction does not occur in some systems, including the RSWEs of section \ref{sec:1drswe}. When $\mathcal{L}^{-1}$ is undefined, a Moore-Penrose pseudoinverse \cite{penrose_inverse} of $\mathcal{L}^{+}$ is used, as suggested in \cite{Tao_simp_imp}. \par
The initial condition for the mean corrected modulation variable in (\ref{eq:aved_meancor_modvar_system}) is defined implicitly. This can be cheaply computed using a fixed-point iteration, given in Algorithm \ref{alg:fixed_point_w0}, starting from a guess of $\vec{W}_0 = \vec{U}_0$. This implicit equation only needs to be solved once, as after the initial condition in the $\vec{W}$ space is found, an explicit timestepping scheme is applied as per standard phase-averaging. \par

\begin{algorithm}
\caption{The Mean Corrected Method}
\label{alg:meancor_alg}
We present the algorithm in terms of Butcher tableau coefficients $a_{ij}, b_i, c_i$ for a general explicit Runge-Kutta scheme \cite{butcher2016numerical}, i.e. classical RK4.
\begin{algorithmic}[1]
    \State Using the standard initial condition, $\vec{U}_0$, compute the initial mean correction, $\vec{C}_0$, and the initial mean corrected modulation variable, $\vec{W}_0$. This usually requires iteration (Algorithm \ref{alg:fixed_point_w0}).
    \For{each timestep}
        \State $\vec{C}_1 \gets \vec{C}_0$
        \State $\vec{W}_1 \gets \vec{W}_0$
        \State Compute $\vec{W}^{n+1}$ from the current time solution, $\vec{W}^{n}$, using phase-averaged function evaluations, $\vec{\overline{f}}_i(\vec{W}_i,t_i)$, for the $M$ intermediate steps i.e.  $M=4$ for RK4.
        \For{$i \in \{ 1, ... M \}$}
            \State $t_i \gets t^n + c_i \Delta t$
            \If{$i \neq 1$}
                \State$\vec{C}_i \gets  \vec{C}(\vec{W}_i, t_i)$ (Use Algorithm \ref{alg:num_C} for a local mean correction)
            \EndIf
            \ParFor{$k \in \{ 1, 2, ..., K \}$ }
                \State $\vec{U}_k \gets  \exp (- ((t_i + s_k) / \epsilon) \mathcal{L}) \vec{W}_i + \epsilon \mathcal{L}^{-1} \vec{C}_i $
                \State $\mathcal{N}_k \gets \mathcal{N}(\vec{U}_k)$
                \State $\vec{f}_k \gets \exp ( ((t_i + s_k) / \epsilon) \mathcal{L})  [\mathcal{N}_k - \mathcal{L} \mathcal{L}^{-1} \vec{C}_i]$
            \EndParFor
            \State $\vec{\overline{f}}_i \gets \sum_{k=1}^{K} \mathcal{K}(s_k) \vec{f}_k$  
            \If{$i \neq M$}
                \State $\vec{W}_{i+1} = \vec{W}^n + \Delta t \sum_{j=1}^{i} a_{i+1,j} \vec{\overline{f}}_i $
            \EndIf
        \EndFor
        \State $\vec{W}^{n+1} = \vec{W}^n + \Delta t \sum_{i=1}^{M} b_i \vec{\overline{f}}_i $ 
        \State$\vec{C}_0 \gets  \vec{C}(\vec{W}^{n+1}, t^{n+1})$
        \State $\vec{W}_0 \gets \vec{W}^{n+1}$
        \State {Back-transform to the standard variable solution at time points of interest:}
        \State $\vec{U}^{n+1} = \exp ( - (t^{n+1}/ \epsilon) \mathcal{L}) \vec{W}^{n+1} + \epsilon \mathcal{L}^{-1} \vec{C}_0  $
    \EndFor
\State Return all values of $\vec{U},\vec{W}$
\end{algorithmic}
\end{algorithm}

\begin{algorithm}
\caption{Fixed-point iteration for computing $\vec{W}(0)$}
\label{alg:fixed_point_w0}
w0\_iteration($\vec{U}_0, C_{\text{tol}}$):
\begin{algorithmic}[1]
    \State$\vec{C}_0 \gets \vec{C}(\vec{U}_0,0)$
    \State$\vec{W}_0 \gets \vec{U}_0 - \epsilon \mathcal{L}^{-1} \vec{C}_0 $ 
    \State$\vec{C}_{\text{new}} \gets \vec{C}(\vec{W}_{0},0) $
    \While{$\sum |\vec{C}_0 - \vec{C}_{\text{new}}| > C_{\text{tol}} $}
        \State$\vec{C}_0 \gets \vec{C}_{\text{new}}$
        \State$\vec{W}_0 \gets \vec{U}_0 - \epsilon \mathcal{L}^{-1} \vec{C}_0 $ 
        \State$\vec{C}_{\text{new}} \gets \vec{C}(\vec{W}_0,0)$
    \EndWhile
\State$\vec{C}_0 \gets \vec{C}_{\text{new}}$
\State Return $\vec{W}_0, \vec{C}_0$
\end{algorithmic}
\end{algorithm}

\subsection{A simple example with the mean correction}
\label{subsec:simple_meancor_example}
We now consider an ODE example with an analytical solution to gain some intuition around the mean corrected version of phase-averaging. Using the idea of \cite{Tao_simp_imp} we consider the simplest `nonlinearity' --- that of a constant, $\mathcal{N} = c$. Our example varies from Tao's in that we consider finite phase-averaging in an oscillatory multiscale system. In this simple ODE, we can compute an analytical solution to compare the accuracy of mean corrected and standard phase-averaging. \par
Consider a standard form (\ref{eq:standard_eq}) ODE for a single variable, with an oscillation frequency $\omega$ and constant $c$: 
\begin{equation}
    \dbydt{u} + \frac{1}{\epsilon} i \omega u = c, ~u(0)=u_0. \label{eq:simp_ex_ODE}
\end{equation}

\noindent We then define the modulation variable transform ($\vec{U} = \exp((t/\epsilon) \mathcal{L}) \vec{V}$),
\begin{equation}
    u = e^{- \frac{i \omega t}{\epsilon}} v,
\label{eq:simp_example_modvar_map}
\end{equation}

\noindent which leads to an initial value problem of
\begin{equation}
    \dbydt{v} = e^{\frac{i \omega t}{\epsilon}} c, ~v(0)=u_0.
\label{eq:example_dv_dt}
\end{equation}

\noindent The analytical solution in the modulation variable space is
\begin{equation}
     v(t) = \frac{i c \epsilon}{\omega} \left( 1 - e^{\frac{i \omega t}{\epsilon}} \right) + u_0. 
     \label{eq:simp_ex_v_sol}
\end{equation}

\noindent Applying the back-transformation (\ref{eq:simp_example_modvar_map}) obtains the standard variable solution:
\begin{equation}
    u(t) = \frac{i c \epsilon}{\omega} \left( e^{-\frac{i \omega t}{\epsilon}} - 1 \right) + e^{-\frac{i \omega t}{\epsilon}} u_0. 
    \label{eq:simp_ex_u_sol}
\end{equation}

Now, consider the application of standard phase-averaging to this ODE:
\begin{equation}
    \dbydt{\overline{v}} = \chi e^{\frac{i \omega t}{\epsilon}} c, ~\overline{v}(0)=u_0, \label{eq:example_dv_bar_dt}
\end{equation}

\noindent with $\chi$ denoting a finite phase-average of the matrix exponential:
\begin{equation}
    \chi = \int_{s=-\eta/2}^{\eta/2} \mathcal{K}(s) e^{\frac{i \omega s}{\epsilon}} ~\text{d}s.
\end{equation}

\noindent When $\eta=0$, $\chi = 1$ and we return to the unaveraged modulation variable equation (\ref{eq:example_dv_dt}). Applying some level of phase-averaging, $\eta=0$, leads to $\chi < 1$ and a reduction in numerical stiffness. Solving the phase-averaged system (\ref{eq:example_dv_bar_dt}) leads to a modulation variable solution of
\begin{equation}
    \overline{v}(t) = \chi \frac{i c \epsilon}{\omega} \left( 1 - e^{\frac{i \omega t}{\epsilon}} \right) + u_0,
\end{equation}

\noindent which in the standard solution space is
\begin{equation}
    \overline{u}(t) = \chi \frac{i c \epsilon}{\omega} \left( e^{-\frac{i \omega t}{\epsilon}} - 1 \right) + e^{-\frac{i \omega t}{\epsilon}} u_0.
\label{eq:simp_example_phase_ave_sol_u}
\end{equation}

\noindent Comparing $\overline{u}(t)$ to the analytical solution of $u(t)$ in (\ref{eq:simp_ex_u_sol}) shows that the averaging error with standard phase-averaging is
\begin{equation}
     u(t) - \overline{u}(t) = (1 - \chi) \frac{i c \epsilon}{\omega} \left( e^{-\frac{i \omega t}{\epsilon}} - 1 \right) .
    \label{eq:simp_example_average_error}
\end{equation}

\noindent This result displays the characteristic trade-off when using finite phase-averaging: an averaging window reduces the numerical stiffness and the associated timestepping error but at the cost of an introduced averaging error. \par
To reduce this averaging error, we turn to the mean corrected method. In this example with $\mathcal{N}=c$, the local mean correction equals the constant $c$ for any choice of mean correction averaging window, $\eta_C$. Noting that $\mathcal{L}^{-1} = (i \omega)^{-1} = -i \omega^{-1}$, we define the mean corrected mapping ($\vec{U} = \exp((-t/\epsilon) \mathcal{L}) \vec{W} + \epsilon \mathcal{L}^{-1} \vec{C}$) of
\begin{equation}
    u = e^{- \frac{i \omega t}{\epsilon}} w - \epsilon \frac{i}{\omega} c,
\end{equation}

\noindent and initial condition,
\begin{equation}
    w_0 = u_0 + \frac{i c \epsilon}{\omega} .
\end{equation}

\noindent The corresponding mean corrected system (\ref{eq:aved_meancor_modvar_system}) is,
\begin{equation}
    \dbydt{\overline{w}} = \left< e^{\frac{i \omega (t+s)}{\epsilon}} [c - c] \right>_s = 0.
\end{equation}

\noindent Hence, the solution in the new modulation variable space is constant in time: $\overline{w}(t) = w_0$. Note, this has not required the approximation of dropping the $\mathcal{O}(\epsilon)$ term in (\ref{eq:meancor_modvar_with_eps_d_dw}) to derive the mean corrected system, as this term is zero regardless: $\epsilon \frac{\partial}{\partial \vec{W}} (\mathcal{L}^{-1} \vec{C}) = \epsilon \frac{\partial}{\partial \vec{W}} (- i \omega^{-1} c) = 0$. \par
With a back-transformation, we obtain the mean corrected solution in the standard space:
\begin{equation}
    \overline{u}(t) =  \frac{i c \epsilon}{\omega} \left( e^{-\frac{i \omega t}{\epsilon}} - 1 \right) + e^{-\frac{i \omega t}{\epsilon}} u_0. 
\end{equation}

\noindent This is precisely the analytical solution of (\ref{eq:simp_ex_u_sol}), so the mean corrected method has not introduced any averaging error. This simple example shows that by including more information in the mean corrected mapping, there may be the potential for smaller averaging errors. The next sections will test whether this improved accuracy holds in more complex differential equations that require numerical simulation.

\section{ODE example: The swinging spring}
\label{sec:ss}
Our first numerical example is the swinging spring ODE. This models a point mass attached to a spring, enabling vertical oscillations in addition to the typical, pendulum, motion. Two distinct physical normal modes, of swinging (with a frequency of $\omega_R$) and elastic springing (with a frequency of $\omega_Z$) exhibit fast, oscillatory, dynamics. A considerable volume of literature has analysed a 2:1 resonance between the swinging and springing modes to mimic a directly resonant interaction of the linear waves, e.g. \cite{swing_spring,swing_spring_multiscale,swing_spring_monodromy}, and with higher-order phase-averaging in  \cite{bauer_colin_spring}. With this direct resonance, the fast modes interact to construct a slow timescale stepwise precession, where the angle of the pendulum position in the $x$-$y$ plane has large shifts in value at certain times \cite{swing_spring}. One full cycle of the nonlinear phase precession occurs on a slower timescale than the swinging and springing modes, corresponding to $\epsilon \approx 0.01$ \cite{bauer_colin_spring}. This presence of multiple timescales intentionally mimics weather and climate models, i.e. the separation of fast inertia-gravity waves and slow Rossby waves of the RSWEs \cite{lynch_swing_spring_resonant_rossby}. \par
In preparation for the PDE systems, which contain a spectrum of in- and out-of-resonance interactions, we examine the swinging spring over a wider range of dynamical scenarios than just a 2:1 resonance. Introducing a `resonance factor' as a ratio of the two normal mode frequencies, $\rho = \omega_Z / \omega_R$, we examine the range of interactions, $\rho \in [1.5,2.5]$. This includes out-of-resonance (e.g. $\rho = 1.5$) and near-resonance (e.g. $\rho = 1.95$) interactions, as well as the 2:1 directly resonant scenario ($\rho = 2$). 

\subsection{Equations}
We use the governing equations in \cite{swing_spring}, which were derived using a cubic-order truncated Lagrangian for the pendulum's position as a function of time,
\begin{subequations}
\begin{align}
    \frac{\text{d}^2 x}{\dt^2} + \omega_R^2 x &= \lambda xz,  \\
    \frac{\text{d}^2 y}{\dt^2} + \omega_R^2 y &= \lambda yz, \\
    \frac{\text{d}^2 z}{\dt^2} + \omega_Z^2 z &=  \frac{\lambda}{2}(x^2 + y^2),
\end{align}
\end{subequations}

\noindent where $\omega_R = \sqrt{g / l}, \omega_Z = \sqrt{k / m}, \lambda = l_o \omega_Z^2 / l^2$, and $l$, $l_0$ are the equilibrium and unstretched spring lengths respectively. We also have the relationship between the two linear frequencies of $\omega_Z = \rho \omega_R$. \par
Following \cite{bauer_colin_spring}, the following transformation is applied to form a first-order system in time:
\begin{equation}
    \dbydt{x} = \omega_R p_x, \
    \dbydt{y} = \omega_R p_y, \
    \dbydt{z} = \omega_Z p_z.
\end{equation}

To work with three variables instead of six, we use a complex-valued solution vector of $\vec{U} = [\hat{x},\hat{y},\hat{z}]^\text{T}$, with
\begin{equation}
    \hat{x} = x + i p_x, \
    \hat{y} = y + i p_y, \
    \hat{z} = z + i p_z.
    \label{eq:complex_sw_spr_defs}
\end{equation}  

\noindent Taking the real part of each component of $\vec{U}$ recovers the position in that coordinate. Using the resonance factor of $\rho =  \omega_Z / \omega_R$ to relate the swinging and springing modes, we arrive at an equation for $\vec{U} = [\hat{x},\hat{y},\hat{z}]^\text{T}$ in the standard form of (\ref{eq:standard_eq}), with linear and nonlinear operators of
\begin{equation}
    \mathcal{L}
    =
    \begin{bmatrix}
i \omega_R & 0 & 0 \\
0 & i \omega_R & 0 \\
0 & 0 & i \rho \omega_R \\
\end{bmatrix}
,
\
\mathcal{L}^{-1}
=
\begin{bmatrix}
\frac{-i}{\omega_R} & 0 & 0 \\
0 & \frac{-i}{\omega_R} & 0 \\
0 & 0 & \frac{-i}{\rho \omega_R} \\
\end{bmatrix}
,
\
\mathcal{N}
=
\begin{bmatrix}
\frac{i \lambda}{\omega_R} \text{Re} (\hat{x}) \text{Re} (\hat{z}) \\
\frac{i \lambda}{\omega_R} \text{Re} (\hat{y}) \text{Re} (\hat{z})  \\
\frac{i \lambda}{2 \rho \omega_R} [\text{Re} (\hat{x})^2 + \text{Re} (\hat{y})^2 ]
\end{bmatrix}.
\label{eq:sw_spring_operators}
\end{equation}

\par
Because $\mathcal{L}$ is a diagonal matrix, the matrix exponential is simply the exponentiation of the nonzero elements,
\begin{equation}
    e^{\pm t \mathcal{L}}
    =
    \begin{bmatrix}
e^{ \pm i \omega_R t} & 0 & 0 \\
0 & e^{ \pm i \omega_R t} & 0 \\
0 & 0 & e^{\pm i \rho \omega_R t} \\
\end{bmatrix} .
\label{eq:sw_spring_e_tL}
\end{equation}

Although $\epsilon$ is not explicitly shown in the swinging spring system (so $\epsilon=1$ is used with Algorithm \ref{alg:meancor_alg}) there is a separation between fast oscillations and slow precession that corresponds to $\epsilon \approx 0.01$ at direct resonance \cite{bauer_colin_spring}.

\subsection{A classical mean correction}
In this first numerical experiment, we consider a classical mean correction corresponding to (\ref{eq:C_inf}). Consider the nonlinear vector in terms of the mean corrected modulation variable, $\vec{W} = [\tilde{x},\tilde{y},\tilde{z}]^\text{T}$,
\begin{equation}
    \mathcal{N}(e^{-t \mathcal{L} }\vec{W})
=
\begin{bmatrix}
\frac{i \lambda}{\omega_R} \text{Re} (e^{-i \omega_R t} \tilde{x}) \text{Re} (e^{-i \rho \omega_R t} \tilde{z}) \\
\frac{i \lambda}{\omega_R} \text{Re} (e^{-i \omega_R t} \tilde{y}) \text{Re} (e^{-i \rho \omega_R t} \tilde{z})  \\
\frac{i \lambda}{2 \rho \omega_R} [\text{Re} (e^{-i \omega_R t} \tilde{x})^2 + \text{Re} (e^{-i \omega_R t} \tilde{y})^2 ]
\end{bmatrix}.
\end{equation}

\noindent Using that $ \text{Re} ( AB ) = (1/2) (AB + A^* B^* )$, where the $^*$ denotes complex conjugation,
\begin{equation}
\mathcal{N}(e^{-t \mathcal{L} }\vec{W})
=
\begin{bmatrix}
\frac{i \lambda}{4 \omega_R} (A \tilde{x} \tilde{z} + B \tilde{x} \tilde{z}^* + C \tilde{x}^* \tilde{z} + D \tilde{x}^* \tilde{z}^*) \\
\frac{i \lambda}{4 \omega_R} (A \tilde{y} \tilde{z} + B \tilde{y} \tilde{z}^* + C \tilde{y}^* \tilde{z} + D \tilde{y}^* \tilde{z}^*)  \\
\frac{i \lambda}{8 \rho \omega_R} (E (\tilde{x}^2 + \tilde{y}^2) + 2(\tilde{x}\tilde{x}^* + \tilde{y}\tilde{y}^*) + F (\tilde{x}^{*2} + \tilde{y}^{*2}) )
\end{bmatrix},
\label{eq:swing_spring_N_w}
\end{equation}

\noindent where,
$$A = e^{-i (\rho + 1) \omega_R t}, \ B = e^{i (\rho - 1) \omega_R t}, \ C = e^{-i (\rho - 1) \omega_R t},$$
$$\ D = e^{i (\rho + 1) \omega_R t}, \ E = e^{-2 i \omega_R t}, \ F = e^{2i \omega_R t}.$$
\par
For the range of resonance factors considered, all the complex exponential terms ($A$ to $F$) will have a nonzero argument. A cancellation of oscillations occurs when computing $\left< . \right>_r$ in the classical limit of $r \rightarrow \infty$, so that each of these terms has a zero mean value (refer to Appendix A for a demonstration of this). As a result, the only nonzero term is in the z-coordinate,
\begin{equation}
    \vec{C}_{\infty} = 
\begin{bmatrix}
0\\
0 \\
\frac{i \lambda}{4 \rho \omega_R} (\tilde{x}\tilde{x}^* + \tilde{y}\tilde{y}^*) 
\end{bmatrix}.
\end{equation}

\noindent This mean correction has a physical interpretation: the positive definite quantity  $\tilde{x} \tilde{x}^* + \tilde{y} \tilde{y}^* = |\tilde{x}|^2 + |\tilde{y}|^2$ is the squared radial position of the pendulum in the modulation variable. \par
As the mean correction only has a component in the $z$ coordinate, the initial states of $\vec{U}(0)$ and $\vec{W}(0)$ have the same $x$ and $y$ components. This allows for the mean corrected initial condition to be defined explicitly from the standard initial condition, $\vec{U}_0 = [\hat{x}_0,\hat{y}_0,\hat{z}_0]^\text{T}$, as
\begin{equation}
    \vec{W}(0)
=
\vec{U}_0 - \mathcal{L}^{-1} \vec{C}_{\infty}(\vec{U}_0)
=
\begin{bmatrix}
 \hat{x}_0 \\
 \hat{y}_0 \\
 \hat{z}_0 - \frac{\lambda}{4 \rho^2 \omega_R^2} (\hat{x}_0 \hat{x}_0 ^* + \hat{y}_0 \hat{y}_0 ^*) \\
\end{bmatrix}.
\end{equation}

\subsection{Numerical implementation}
Parameters following \cite{lynch_swing_spring_resonant_rossby} of $g = \pi^2 ~\text{m} ~\text{s}^{-2}, ~l = 1 ~\text{m}, ~m = 1 ~\text{kg}$ are used, which means that $\omega_R = \pi ~\text{s}^{-1}$. We additionally choose that $l_0 = 1.2 ~\text{m}, k = \rho^2 \pi^2 ~\text{kg} ~\text{s}^{-2}$, such that $\omega_Z = \rho \omega_R$. The initial condition is $\vec{U}(0) = [0.04, 0.03427 i / \pi,0.08]$ and the simulations are run for a total time of 200 s. Explicit RK4 is used to solve (\ref{eq:modvar_angle_bracket_ave}) and (\ref{eq:meancor_modvar_angle_bracket_ave}) for the phase-averaged and mean corrected (with Algorithm \ref{alg:meancor_alg}) methods respectively, using the operators in (\ref{eq:sw_spring_operators}). The large timestep solutions use $\Delta t = 0.5 ~\text{s}$. Errors are measured as an average L2 difference, over the simulation, in the coordinate position relative to a reference solution computed using RK4 and $\Delta t = 0.01 ~\text{s}$. \par

\subsection{Results}
\begin{figure}
    \centering
     \includegraphics[scale = 0.7]{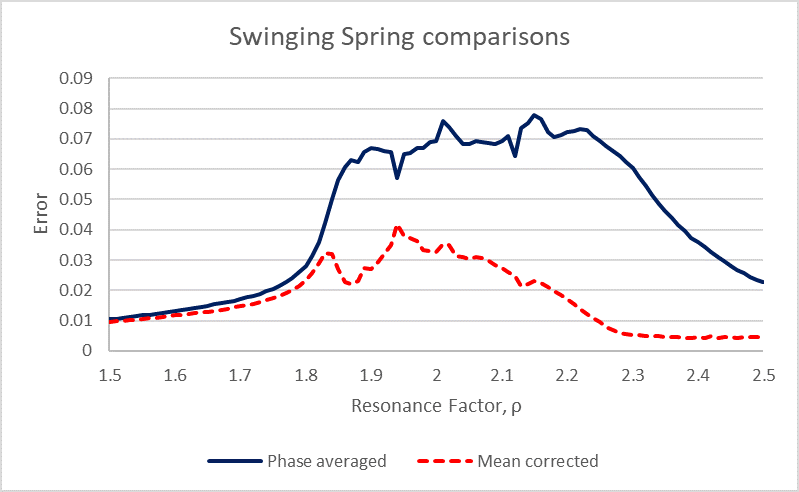}
    \caption{Errors for standard (solid line) and mean corrected (dashed) phase-averaging in the swinging spring experiment. The mean corrected method is more accurate for all the resonance factors examined, which models a range of in- and out-of-resonance interactions. }
    \label{fig:sw_spring_results}
\end{figure}

Errors at the best averaging window for each method are shown in Figure \ref{fig:sw_spring_results}. The mean corrected method has a lower error for all resonance factors, so improves accuracy in the directly resonant, near-resonant and out-of-resonance interaction settings of the swinging spring. This efficacy for a range of dynamical situations motivates the application of a mean correction in PDEs with a bilinear $\mathcal{N}(\vec{U},\vec{U})$ like the cases in sections \ref{sec:kg} and \ref{sec:1drswe}, as these systems can be expanded as a sum of three-wave interactions of different frequencies \cite{Embid_1996, Adam}. The accuracy improvement is greater for interactions of $\rho > 1.82$, which highlights that the mean correction may have a larger impact for certain regimes of nonlinear dynamics.

\par 
\par
\begin{figure}
     \centering
     \begin{subfigure}[b]{0.48\textwidth}
         \centering
         \includegraphics[width=\textwidth]{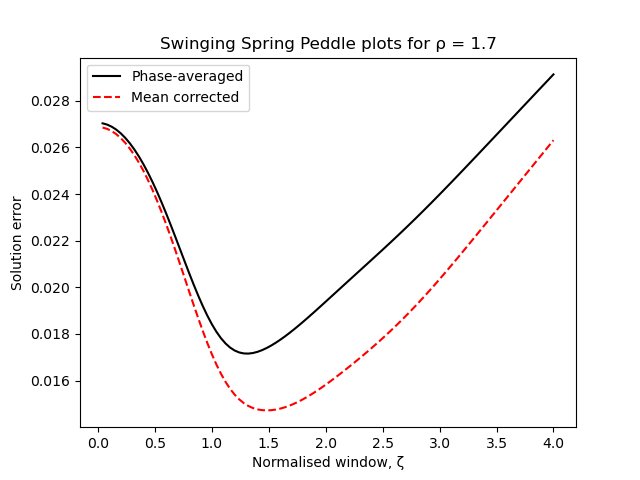}
         \caption{Out of resonance ($\rho = 1.7$)}
     \end{subfigure}
     \hfill
     \begin{subfigure}[b]{0.48\textwidth}
         \centering
         \includegraphics[width=\textwidth]{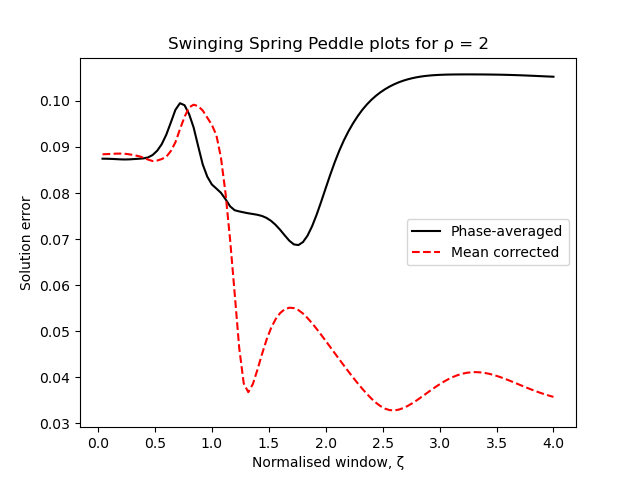}
         \caption{In resonance ($\rho = 2$)}
         \label{fig:swing_spring_peddle_rho_2}
     \end{subfigure}
        \caption{Peddle plots of the errors over normalised averaging window, $\zeta = \eta /\Delta t$. In an out-of-resonance case (a), there is a clear minimum in the error for each method, with the best window for the mean corrected method being slightly larger. For the directly resonant case (b), the mean corrected method again is best with a larger window, but there are also multiple local minima in the error.}
    \label{fig:swing_spring_peddle}
\end{figure}

The best averaging windows are identified from the lowest error of each method over a range of averaging windows separated by $\Delta \zeta = 0.1$ ($\zeta = \eta / \Delta t$). The existence of optimal averaging windows is shown in Peddle plots in Figure \ref{fig:swing_spring_peddle}. The mean corrected method often has multiple minima in these plots, e.g. Figure \ref{fig:swing_spring_peddle_rho_2}. 

\par

\par
\begin{figure}
     \centering
     \begin{subfigure}[b]{0.47\textwidth}
         \centering
         \includegraphics[width=\textwidth]{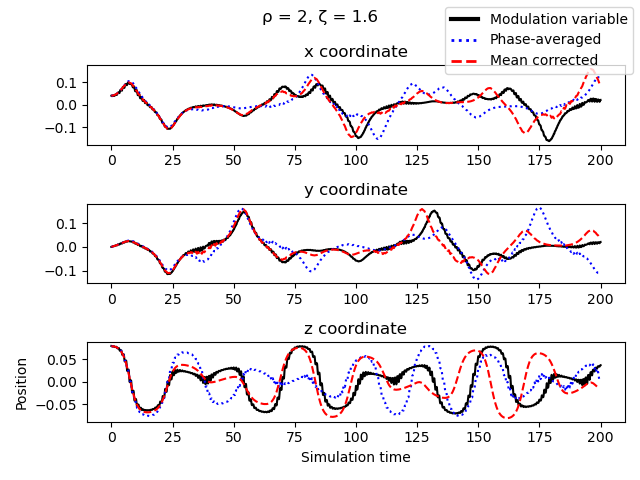}
         \caption{Best $\zeta$ for the phase-averaged method }
         \label{fig:swing_spring_v_best}
     \end{subfigure}
     \hfill
     \begin{subfigure}[b]{0.47\textwidth}
         \centering
         \includegraphics[width=\textwidth]{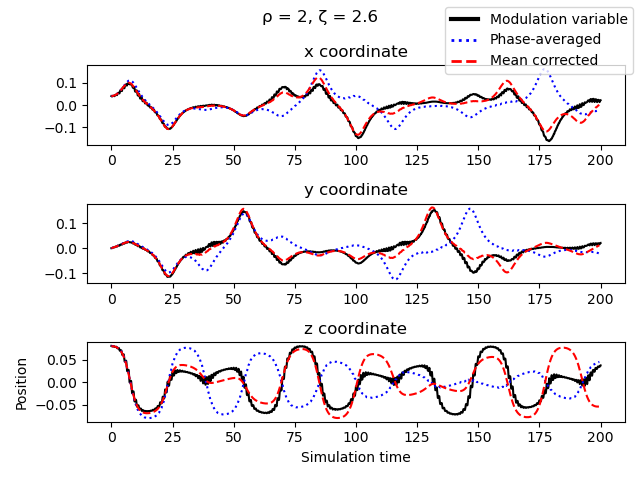}
         \caption{Best $\zeta$ for the mean corrected method}
         \label{fig:swing_spring_w_best}
     \end{subfigure}
        \caption{Comparing the unaveraged modulation variable solution ($\vec{V}$, solid line), with the modulation variables of the standard phase-averaged method ($\vec{\overline{V}}$, dotted) and mean corrected method ($\vec{\overline{W}}$, dashed). This is in the directly resonant state of $\rho = 2$. The modulation variable solutions in (a) use the best $\zeta$ for the phase-averaged method, whilst (b) use the best $\zeta$ for the mean corrected method. The mean corrected method aligns better in the modulation space with its best averaging window in (b), whereas the phase-averaged method has considerable phase discrepancies, even in (a).}
        \label{fig:swing_spring_modvar_space}
\end{figure}

The effect of the mean correction can be observed when comparing modulation variable spaces in Figure \ref{fig:swing_spring_modvar_space}. An unaveraged modulation variable solution using (\ref{eq:modvar_system}) represents the `true' modulation variable solution. This still contains small amplitude fluctuations, which the phase-averaged method seeks to reduce. This introduces an averaging error, which presents itself as a phase discrepancy in the modulation space, even at the best window for standard phase-averaging (Figure \ref{fig:swing_spring_v_best}). With a larger averaging window (Figure \ref{fig:swing_spring_w_best}), the standard phase-averaged solution is smoother but the phase errors are more prevalent, with the z-coordinate solution often being 180 degrees out of phase. By contrast, the mean corrected method has a much lower phase error with its best $\eta^*$ in Figure \ref{fig:swing_spring_w_best}. Hence, the information provided by the mean correction increases the accuracy within the modulation variable space itself. 

\section{PDE system 1: One-dimensional Klein-Gordon-type PDE}
\label{sec:kg}
\subsection{Equations}
We now turn to a PDE example of a one-dimensional Klein-Gordon-type (KG-type) equation. This is one of the simplest wave systems with purely oscillatory linear dynamics. We consider a quadratic nonlinearity,
\begin{equation}
    \frac{\partial^2 u}{\partial t^2} - \alpha \frac{\partial^2 u}{\partial x^2} + \beta u + \gamma u^2 = 0.
\end{equation}

\noindent A multiscale system can be generated with a specific choice of parameters $\alpha, \beta, \gamma$,
\begin{equation}
    \frac{\partial^2 u}{\partial t^2} + \frac{1}{\epsilon^2} \left( u - \frac{\partial^2 u}{\partial x^2} \right) = - u^2. 
    \label{eq:kg_pde}
\end{equation}

A linear dispersion relation of $\omega = \pm \sqrt{1 + k^2}$ can be identified from (\ref{eq:kg_pde}) in Fourier space, with the right-hand side set to zero and $\epsilon = 1$. To obtain a first-order-in-time system with a skew-Hermitian $\mathcal{L}$, we introduce variables of $a= ( \omega / \epsilon ) u, ~b = \frac{\partial u }{ \partial t}$. For a solution vector of $\vec{U} = [a,b]^\text{T}$, we have the following operators in Fourier space for a system in the standard form (\ref{eq:standard_eq}),
\begin{equation}
\mathcal{L} 
= 
\begin{bmatrix}
0 & -\omega\\
\omega & 0\\
\end{bmatrix}
,
\
\mathcal{L}^{-1}
=
\begin{bmatrix}
0 && \frac{1}{\omega} \\
\frac{-1}{\omega} && 0
\end{bmatrix}
,
\
\mathcal{N}  
= 
\begin{bmatrix}
0\\
- \left( \frac{\epsilon}{\omega} a \right) \circledast \left( \frac{\epsilon}{\omega} a \right) 
\end{bmatrix},
\label{eq:KG_operators}
\end{equation}
where $\circledast$ denotes a circular spatial convolution. In this standard form, the KG-type equation is both weakly nonlinear ($\mathcal{N} \sim \mathcal{O}(\epsilon^2)$) and multiscale ($\epsilon^{-1} \mathcal{L} \sim \mathcal{O}(\epsilon^{-1})$). \par
The eigenvalues of the linear operator are $\pm i \omega$, with the corresponding Fourier space matrix exponential only containing trigonometric terms,
\begin{equation}
    e^{\pm \frac{t \mathcal{L}}{\epsilon}}
    =
    \begin{bmatrix}
    \cos(\frac{\omega t}{\epsilon}) & \mp \sin(\frac{\omega t}{\epsilon})\\
    \pm \sin(\frac{\omega t}{\epsilon}) & \cos(\frac{\omega t}{\epsilon})\\
    \end{bmatrix}.
    \label{eq:KG_matexp}
\end{equation}

We consider a classical mean correction (\ref{eq:C_inf}), like was used for weakly nonlinear systems in \cite{Tao_simp_imp}, as well as a local mean correction (\ref{eq:C_finite}). Details of the classical mean correction, which is computed in Fourier space, are given in Appendix B. The local mean correction is computed using Algorithm \ref{alg:meancor_alg} with the computationally cheap choice of $\eta_C = \eta$. 

\subsection{Numerical implementation}
A one-dimensional spatial domain of length $[0,2 \pi]$ is used with $N_x = 32$ grid points and periodic boundary conditions. A Gaussian perturbation is applied to the $a$ field:
\begin{equation}
    a(x,0) =  \exp \left( \frac{-(x-\pi)^2}{2} \right), ~b(x,0) = 0.
\label{eq:KG_IC}
\end{equation}

We use a pseudospectral method \cite{boyd2001chebyshev,canuto2007spectral} for this system and the RSWEs in section 5. This method computes derivatives in spectral space and multiplication in physical space, with the use of Fast Fourier Transforms to map between spaces. Advantageously, this spatial discretisation allows Fourier space operators (\ref{eq:KG_operators}) and an analytical matrix exponential (\ref{eq:KG_matexp}) to be used, as well as the classical mean correction in Appendix B. \par
Using the operators in (\ref{eq:KG_operators}), an explicit RK4 scheme is applied to equation (\ref{eq:modvar_angle_bracket_ave}) for the phase-averaged method, with the number of kernel points being defined using (\ref{eq:K_rule}) and $\omega_{\text{max}} = \sqrt{1+k_{\text{max}}^2}$. Algorithm \ref{alg:meancor_alg} is used for the mean corrected method. A total simulation time of $T_{\text{\text{max}}} = 20$ is used, with large timesteps of $\Delta t \in [1,3]$. The error metric is the average L1 difference in the $a$ field over the simulation, relative to an unaveraged modulation variable RK4 solution (\ref{eq:modvar_system}) with $\Delta t = 10^{-4}$. \par
The initial condition for the mean corrected modulation variable is computed using Algorithm \ref{alg:fixed_point_w0} and $C_{\text{tol}} = 10^{-10}$. For these experiments,
$\epsilon \in \{ 0.5, 0.1, 0.05, 0.01  \}$ are examined. We compute the errors of both methods over a range of averaging windows separated by $\Delta \zeta = 0.05$. Different $\eta^*$ are usually identified for all three methods, although these are identical for $\epsilon = 0.01$ for all but the largest timestep. \par

\subsection{Results}
\begin{figure}
        \centering
        \begin{subfigure}[b]{0.475\textwidth}
            \centering
            \includegraphics[width=\textwidth]{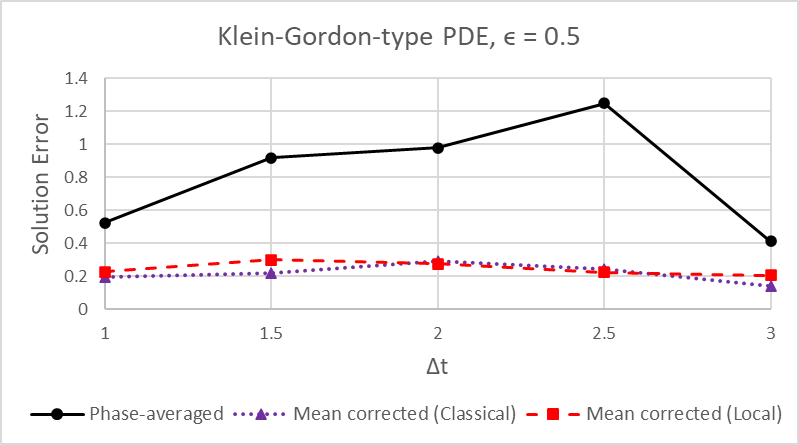}
            \caption{$\epsilon = 0.5$}
        \end{subfigure}
        \hfill
        \begin{subfigure}[b]{0.475\textwidth}  
            \centering 
            \includegraphics[width=\textwidth]{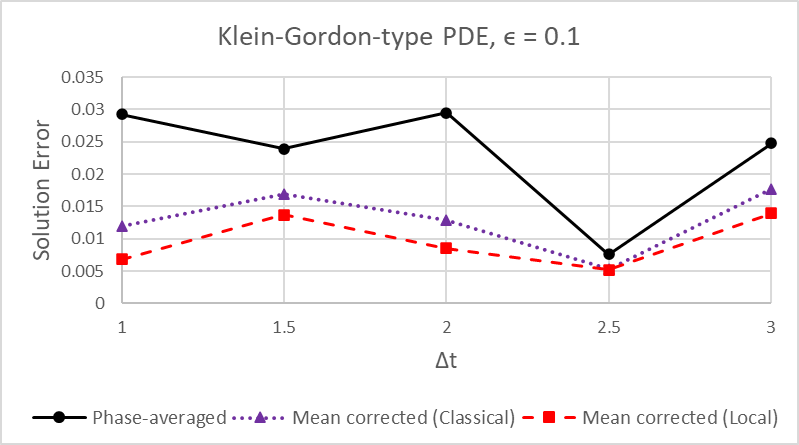}
            \caption{$\epsilon = 0.1$}
        \end{subfigure}
        \vskip\baselineskip
        \begin{subfigure}[b]{0.475\textwidth}   
            \centering 
            \includegraphics[width=\textwidth]{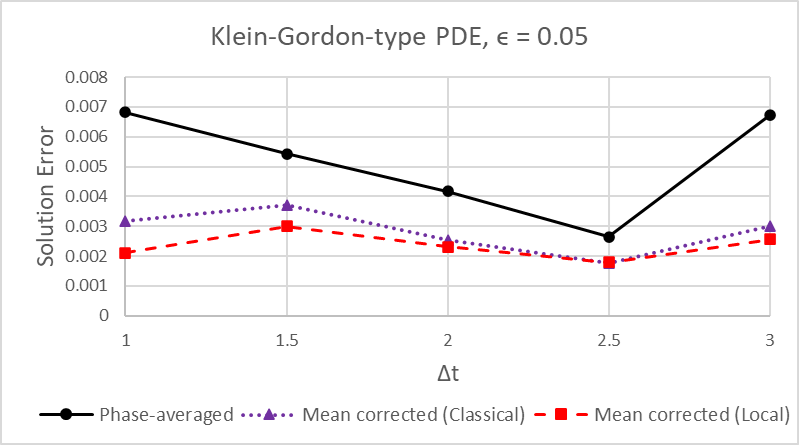}
            \caption{$\epsilon = 0.05$}
        \end{subfigure}
        \hfill
        \begin{subfigure}[b]{0.475\textwidth}   
            \centering 
            \includegraphics[width=\textwidth]{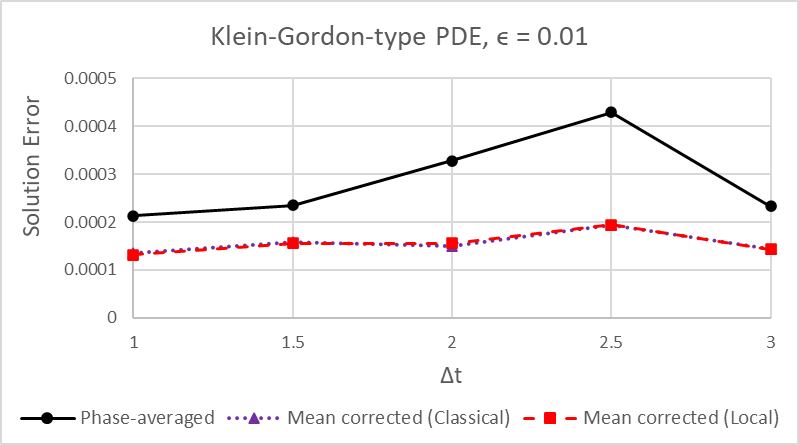}
            \caption{$\epsilon = 0.01$}
        \end{subfigure}
        \caption{Errors for the standard phase-averaged (solid line), classically mean corrected (dotted line), and locally mean corrected (dashed line) methods in the KG-type system. The two mean corrected methods are more accurate than the phase-averaged one for all values of $\epsilon$. The local mean correction has a lower error than the classical for $\epsilon \in \{ 0.1, 0.05 \}$. At $\epsilon=0.01$ there is minimal difference between $\vec{C}_{\infty}$ and $\vec{C}$, so both mean corrected methods perform similarly.}
    \label{fig:KG_results}
\end{figure}

\begin{figure}[ht]
     \centering
     \begin{subfigure}[b]{0.48\textwidth}
         \centering
         \includegraphics[width=\textwidth]{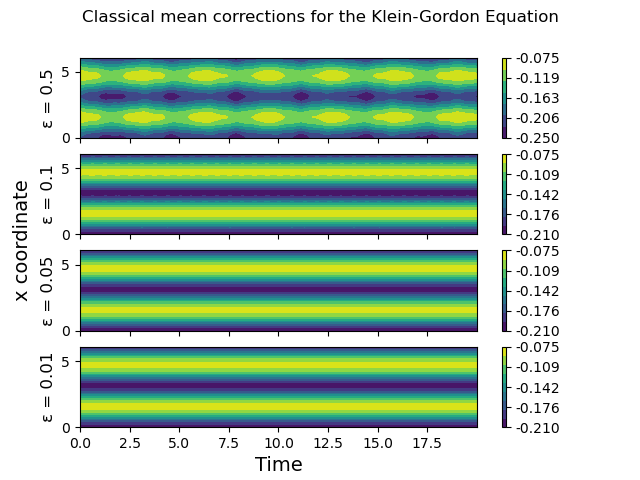}
         \caption{Classical mean corrections}
     \end{subfigure}
     \hfill
     \begin{subfigure}[b]{0.48\textwidth}
         \centering
         \includegraphics[width=\textwidth]{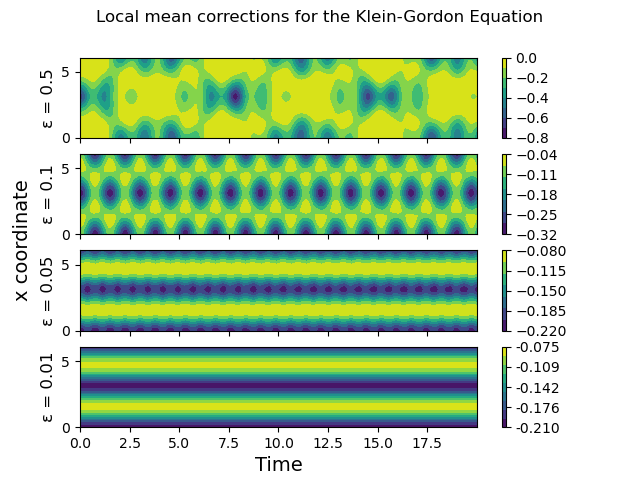}
         \caption{Local mean corrections}
     \end{subfigure}
        \caption{Visualisations of the classical and local mean corrections in the KG-type system. The best $\eta_C$ for each method at $\Delta t = 1$ is used to compute mean corrections from a small timestep unaveraged modulation variable ($\vec{V}$) solution. The inverse Fourier transform is performed to give a physical (as opposed to spectral) space visualisation. Mean corrections from faster linear oscillations are shown moving down the page, with $\epsilon \in \{ 0.5, 0.1, 0.05, 0.01 \}$. The classical mean correction has a different profile for $\epsilon = 0.5$ but is visually indistinguishable in the other regimes.  By contrast, the local mean correction has an observable variation between each $\epsilon$ regime. This is a result of the finite $\eta_C$ capturing different local mean dynamics depending on the oscillation frequency. At $\epsilon = 0.01$, the mean corrections are visually indistinguishable; this coincides with the two methods having similar accuracy.}
    \label{fig:KG_an_num_Cs}
\end{figure}

Figure \ref{fig:KG_results} shows that both implementations of the mean correction (classical and local) have lower errors than the standard phase-averaged method, for the examined timesteps and values of $\epsilon$. The local mean correction has a lower error than the classical for $\epsilon \in \{ 0.1, 0.05 \}$. For the fastest oscillation case of $\epsilon = 0.01$, the mean corrected method errors differ on average by 0.3$\%$ relative to standard phase-averaging. This is due to the local mean correction tending to the classical as $\epsilon$ reduces; we observe this in visualisations of the mean corrections in physical space in Figure \ref{fig:KG_an_num_Cs}. The classical mean correction has a near identical mean correction for all regimes except $\epsilon = 0.5$, whilst the local corrections are visually distinct at each $\epsilon$. This reflects how using a finite averaging window for the mean correction captures local variations in the mean of the oscillations. The classical correction loses this with the complete removal of oscillations. The largest $\epsilon$ value of $0.5$ sees the two methods performing similarly, with the classical correction sometimes being more accurate. A possible reason for this may the neglection of the $\epsilon \frac{\partial}{\partial \vec{W}} (\mathcal{L}^{-1} \vec{C})$ term from (\ref{eq:meancor_modvar_with_eps_d_dw}), as the local correction has a larger magnitude and more variation (Figure \ref{fig:KG_an_num_Cs}) which increases the size of this dropped term. 
\par

\section{PDE system 2: One-dimensional rotating shallow water equations}
\label{sec:1drswe}
\subsection{Equations}
We lastly demonstrate the mean corrected method in the one-dimensional rotating shallow water equations (RSWEs), which are a simplified system relevant for weather and climate applications. These equations were used to test other phase-averaged implementations in \cite{Haut_Wingate_2014,peddle_haut_wingate_2019,rosemeier2024multilevel}. This test will only consider the local version of the mean correction. This is
the more practical option for use in the complex PDEs governing the weather and
climate as it can be computed numerically for any spatial discretisation choice. \par
The prognostic variables of the RSWEs are a two-dimensional fluid velocity, $(u,v)$, and a measure of elevation of the fluid layer, such as a height perturbation, $\eta_h$, above the mean fluid height of $H_0$. We will use the nondimensional $f$-plane RSWEs in one spatial dimension, derived in detail in \cite{Adam}. These use an elevation measure of a symmetrical geopotential, $\phi = \sqrt{g / H_0} ~\eta_h$, where $g$ is the gravitational acceleration. Representing this system in Fourier space, where $\hat{\vec{U}} = [\hat{u},\hat{v},\hat{\phi}]^\text{T} $ contains the spatially Fourier transformed prognostic variables, we have the following operators in the form of our standard equation (\ref{eq:standard_eq}),
\begin{equation}
    \mathcal{L}
    =
    \begin{bmatrix}
0 & -1 & i k\\
1 & 0 & 0\\
i k & 0 & 0\\
\end{bmatrix}
,
\
\mathcal{N}
=
\mathcal{F} 
\Bigg\{ -
\begin{bmatrix}
u\frac{\partial u}{\partial x} \\
u \frac{\partial v}{\partial x}  \\
\frac{\partial }{\partial x} (u \phi) 
\end{bmatrix}
\Bigg\},
\label{eq:rswe_operators}
\end{equation}

\noindent where $\mathcal{F}$ denotes the spatial Fourier transform. \par
The linear dispersion relation of $\mathcal{L}$ in (\ref{eq:rswe_operators}) is $\omega = \pm \alpha \sqrt{1 + k^2 }$, with $\alpha \in \{ -1,0,1 \}$ denoting the dispersion relation branches. The fast modes of $\alpha = \pm 1$ correspond to inertia-gravity waves with the same frequencies as the waves in the KG-type equation (\ref{eq:kg_pde}). Additionally, there is a slow mode of $\omega = 0$ corresponding to Rossby waves. The speed of these Rossby waves is identically zero as a result of a constant Coriolis parameter. A further consequence is that $\mathcal{L}$ is singular and $\mathcal{L}^{-1}$ is undefined. As such, the Moore-Penrose pseudoinverse, $\mathcal{L}^{+}$, is used in place of the regular inverse. This means that $ \mathcal{L} \mathcal{L}^{+} \neq I$,
\begin{equation}
\mathcal{L^+}
=
\begin{bmatrix}
0 & \frac{1}{\psi^2} & \frac{-ik}{\psi^2}\\
\frac{-1}{\psi^2} & 0 & 0\\
\frac{-ik}{\psi^2} & 0 & 0\\
\end{bmatrix}
, 
\
\mathcal{L} \mathcal{L^+}
=
\begin{bmatrix}
1 & 0 & 0\\
0 & \frac{1}{\psi^2} & \frac{-ik}{\psi^2}\\
0 & \frac{ik}{\psi^2} & \frac{k^2}{\psi^2}\\
\end{bmatrix},
\label{eq:rswe_L_Linv}
\end{equation}
\noindent where $\psi = \sqrt{1 + k^2}$. 

The matrix exponential for this system can be written analytically in Fourier space,
\begin{equation}
    e^{\pm \frac{t}{\epsilon} \mathcal{L}}
    =
    \begin{bmatrix}
    \cos(\frac{\psi t}{\epsilon}) & \mp \frac{1}{\psi} \sin(\frac{\psi t}{\epsilon}) & \pm \frac{ik}{\psi} \sin(\frac{\psi t}{\epsilon})\\
    \pm \frac{1}{\psi} \sin(\frac{\psi t}{\epsilon}) & \frac{1}{\psi^2}(k^2  + \cos(\frac{\psi t}{\epsilon})) & \frac{ik}{\psi^2}(1-\cos(\frac{\psi t}{\epsilon}))\\
    \pm \frac{ik}{\psi} \sin(\frac{\psi t}{\epsilon}) & \frac{ik}{\psi^2} (\cos(\frac{\psi t}{\epsilon}) -1) & \frac{1}{\psi^2} (1 + k^2  \cos(\frac{\psi t}{\epsilon}))\\
    \end{bmatrix}.
    \label{eq:matrix_exp_RSWEs}
\end{equation}

Numerical dissipation is required for stabilising the pseudospectral method and is applied in the form of a $(\nabla^2)^2$ hyperviscosity,
\begin{equation}
\mathcal{D} = 
\begin{bmatrix}
-\mu k^4 & 0 & 0\\
0 & -\mu k^4 & 0 \\
0 & 0 & -\mu k^4 \\
\end{bmatrix}.
\end{equation}

We use the following equations for the mean corrected method,
\begin{subequations}
\begin{align}
    \vec{U} = e^{-\frac{t}{\epsilon}\mathcal{L}} \vec{W} + \epsilon \mathcal{L^+} \vec{C}, \\ 
    \frac{\partial \vec{\overline{W}}}{\partial t} = \left< e^{ \frac{t+s}{\epsilon} \mathcal{L}} \left[ \mathcal{N}^* \left( e^{- \frac{t+s}{\epsilon} \mathcal{L}} \vec{\overline{W}}(t) + \epsilon \mathcal{L}^{+} \vec{C} \right) - \mathcal{L} \mathcal{L}^{+} \vec{C}  \right] \right>_s, \\
    \vec{\overline{W}}(0) = \vec{U_0} - \epsilon \mathcal{L}^{+} \vec{C}(\vec{\overline{W}}(0)),
\end{align}
\label{eq:RSWEs_meancor_eq}
\end{subequations}

\noindent where $\mathcal{N}^* (\vec{U}) = \mathcal{N}(\vec{U}) + \mathcal{D} \vec{U}$ is the nonlinear vector with added dissipation. 

\subsection{Numerical implementation}
We consider initial conditions of a Gaussian perturbation to the geopotential height field, at a state of rest,
\begin{equation}
    u(x,0) = 0, \ \ \  v(x,0) = 0, \ \ \ \phi(x,0) = \exp \left( \frac{-(x-\pi)^2}{2} \right).
    \label{eq:rswe_gaussian_IC}
\end{equation}

A spatially periodic domain of size $[0,2 \pi]$ is used with $N_x = 32$ grid points. The initial condition for the modulation variable is computed with the fixed-point iteration of Algorithm \ref{alg:fixed_point_w0} and $C_{\text{tol}} = 10^{-10}$. A hyperviscosity coefficient of $\mu = 10^{-4}$ is used. The phase-averaged solution is computed using (\ref{eq:modvar_angle_bracket_ave}) with the dissipation-stabilised $\mathcal{N}^*$ in place of $\mathcal{N}$. The mean corrected method applies Algorithm \ref{alg:meancor_alg} to (\ref{eq:RSWEs_meancor_eq}). The use of a pseudospectral method allows the analytical matrix exponential (\ref{eq:matrix_exp_RSWEs}) to be used. Alternate matrix exponential constructions for the RSWEs include rational approximations, e.g. \cite{Haut_Babb,schreiber2018beyond}, or Chebyshev approximation \cite{hiroe_colin}. We present results for a total simulation time of $T_{\text{\text{max}}} = 10$. The error metric is the average over the simulation of an L2 difference in solution components relative to a reference RK4 solution with $\Delta t = 10^{-4}$. \par
The locally mean corrected method requires the selection of an ($\eta$, $\eta_C$) pairing. Our ($\eta^*$, $\eta_C^*$) selection procedure for this test is:

\hfill

\begin{enumerate}[i)]
    \item Run simulations with the standard phase-averaged method over a range of $\zeta = \eta / \Delta t$ values to find the window with lowest error, $\eta^*$. For the RSWEs, increments of $\Delta \zeta = 0.05$ are used.
    \item Using the $\eta^*$ value found in i), run simulations with the mean corrected method over a range of $\eta_C$ values to identify the local mean correction window with the lowest error, $\eta_C^*$. For the RSWEs, increments of $ \Delta \eta_C = \epsilon$ are used.
\end{enumerate}

As this process applies the same window of $\eta^*$ to the timestepping equation in both the standard and mean corrected methods, the mean corrected method only needs to be examined over a range of $\eta_C$ values, instead of ($\eta$, $\eta_C$) combinations.

\hfill

\subsection{Results}

\begin{figure}
        \centering
        \begin{subfigure}[b]{0.475\textwidth}
            \centering
            \includegraphics[width=\textwidth]{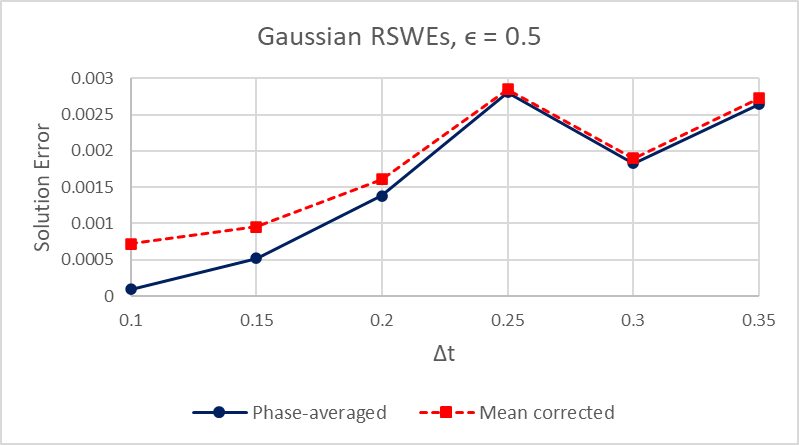}
            \caption{$\epsilon = 0.5$}
        \end{subfigure}
        \hfill
        \begin{subfigure}[b]{0.475\textwidth}  
            \centering 
            \includegraphics[width=\textwidth]{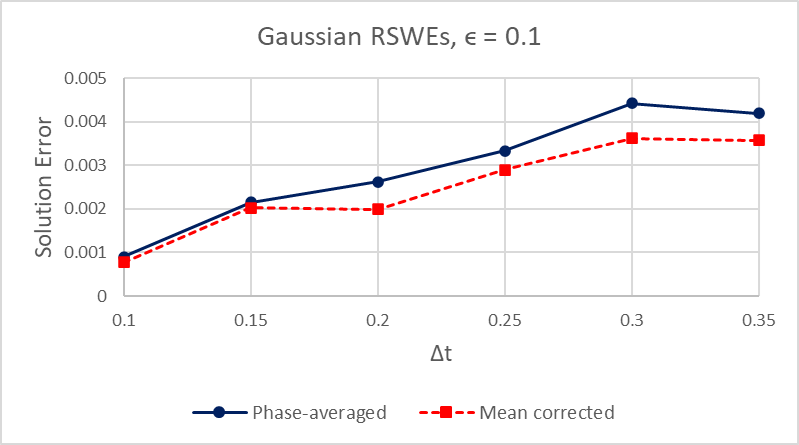}
            \caption{$\epsilon = 0.1$}
        \end{subfigure}
        \vskip\baselineskip
        \begin{subfigure}[b]{0.475\textwidth}   
            \centering 
            \includegraphics[width=\textwidth]{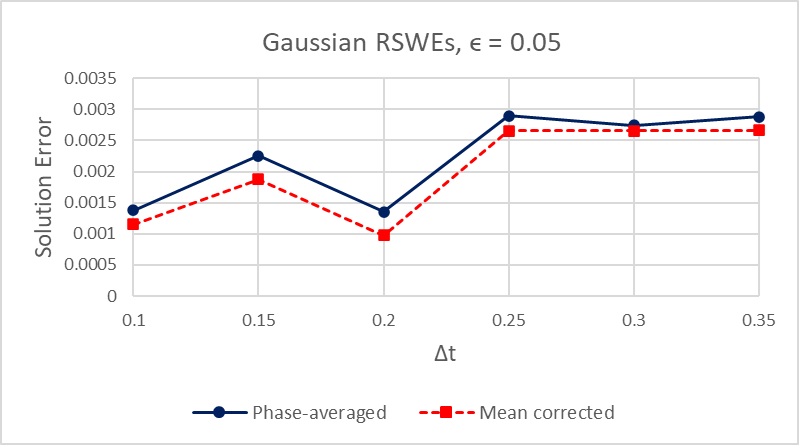}
            \caption{$\epsilon = 0.05$}
        \end{subfigure}
        \hfill
        \begin{subfigure}[b]{0.475\textwidth}  
        
            \centering 
            \includegraphics[width=\textwidth]{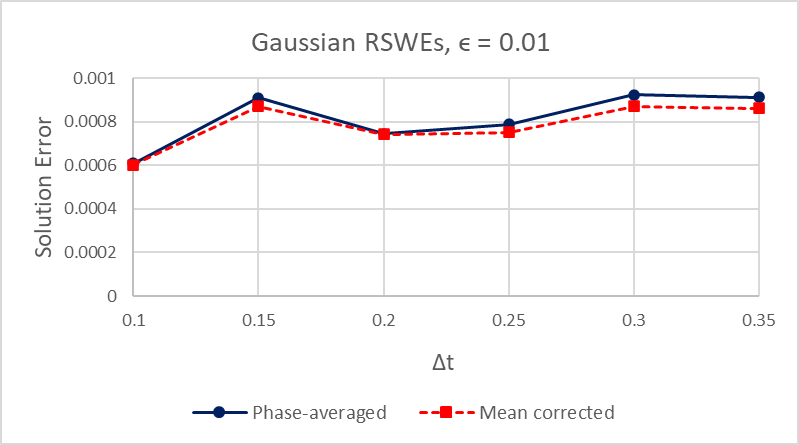}
            \caption{$\epsilon = 0.01$}
        \end{subfigure}
        \vskip\baselineskip
        \begin{subfigure}[b]{0.475\textwidth}   
            \centering 
            \includegraphics[width=\textwidth]{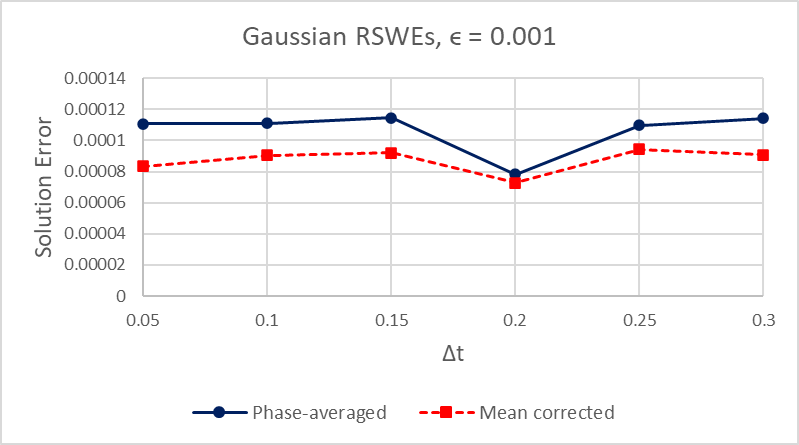}
            \caption{$\epsilon = 0.001$}
        \end{subfigure}
        
        \caption{Standard (solid line) and locally mean corrected (dashed) phase-averaging errors for the RSWEs experiment. Accuracy improvements with the mean correction are observed in the regimes with $\epsilon \leq 0.1$. For this implementation in the RSWEs, the linear oscillations with $\epsilon = 0.5$ are insufficiently fast for the mean correction to be effective. }
    \label{fig:RSWE_results}
\end{figure}

The phase-averaged and locally mean corrected methods errors are compared at timescale separations of $\epsilon \in \{ 0.5, 0.1, 0.05, 0.01, 0.001 \} $ in Figure \ref{fig:RSWE_results}. The local mean correction enables greater accuracy in the cases with sufficiently fast oscillations ($\epsilon  \leq 0.1 $). Although for some tests the mean correction has a minimal effect (e.g. $\Delta t = 0.1, 0.2$ for $\epsilon = 0.01$) it crucially does no worse than standard phase-averaging for these configurations. This is consistent with this proven property for the weakly nonlinear, classically time-averaged case in \cite{Tao_simp_imp}. The mean correction continues to provide accuracy improvements at very large timescale separations ($\epsilon = 0.001$). The mean corrected method does not perform better for the timescale separation of $\epsilon = 0.5$ and a possibility for this is neglecting the $ \mathcal{O} (\epsilon)$ term in (\ref{eq:meancor_modvar_with_eps_d_dw}).  
\par
\par

\begin{figure}
     \centering
     \begin{subfigure}[t]{0.475\textwidth}
         \centering
         \includegraphics[width=\textwidth]{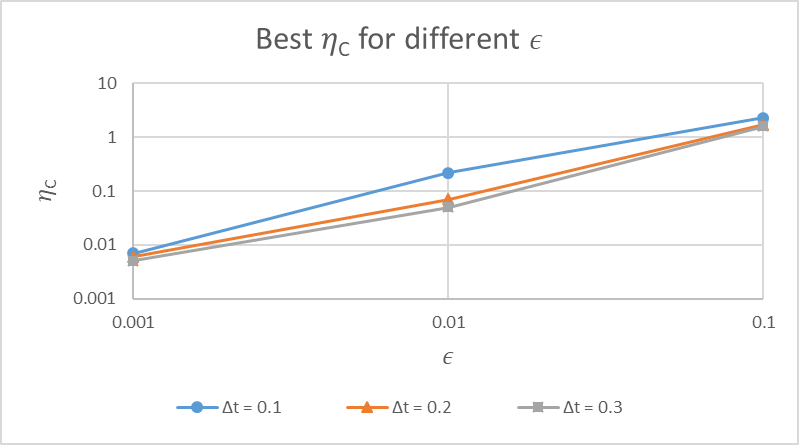}
         \caption{Best $\eta_C$ for different $\epsilon$}
         \label{fig:RSWEs_etaC_over_eps}
     \end{subfigure}
     \hfill
     \begin{subfigure}[t]{0.475\textwidth}
         \centering
         \includegraphics[width=\textwidth]{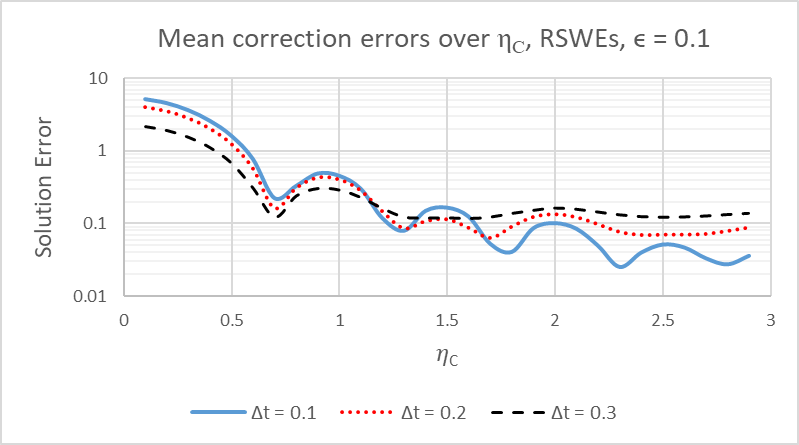}
         \caption{Mean corrected errors over $\eta_C$, with $\epsilon = 0.1$}
         \label{fig:RSWEs_meancor_err_over_etaC}
     \end{subfigure}
        \caption{Investigating the selection of the best  $\eta_C$ for the RSWEs. In (a), the mean correction window with the lowest error, $\eta_C^*$, is plotted against timescale separations of $\epsilon \in \{0.1, 0.01, 0.001 \}$ for different timestep sizes. There is a relatively linear scaling between $\eta_C^*$ and $\epsilon$ on these logarithmic axes. This reflects that the size of $\eta_C^*$ adjusts to consider a similar number of linear oscillations within the support of $\mathcal{K}_C$. In (b), the local minimal error locations are similar with different timestep sizes, although a larger $\Delta t$ leads to less variation in the error over $\eta_C$. The best $\eta_C$ generally reduces with larger timestep: the third local minima around $\eta_C = 2.7$ is best for $\Delta t = 0.1$, whilst the second local minima of $\eta_C = 1.7$ is best for $\Delta t \in \{ 0.2, 0.3 \}$.}
    \label{fig:RSWE_etaC_error_plots}
\end{figure}

In Figure \ref{fig:RSWE_etaC_error_plots} we examine the best value of $\eta_C$ in these experiments. With our use of a finite level of phase-averaging ($\eta^*$), the mean correction window that leads to the lowest error ($\eta_C^*$) is also one of finite length. This shows the potential for the local mean correction, which retains low-frequency nonlinear oscillations, to be more effective than the classical version for weather and climate applications. \par
The size of $\eta_C^*$ is dependent on the timescale separation (Figure \ref{fig:RSWEs_etaC_over_eps}). As $\epsilon$ reduces, a smaller averaging window is needed to capture a similar number of linear oscillations in the support of $\mathcal{K}_C$. As such, the computation of the mean correction scales well to larger timescale separations, as the required $K_r$ for the $\left< \cdot \right>_r$ operation does not increase with a smaller $\epsilon$, whilst $K_s$ does for the standard $ \left< \cdot \right>_s $ phase-average (\ref{eq:K_rule}). The errors for the mean corrected method over $\eta_C$ differ from those seen in the PDE Peddle plots over $\eta$, in that there are multiple local minima as opposed to a single optimum (Figure \ref{fig:RSWEs_meancor_err_over_etaC}). The locations of these minima are similar for different values of $\Delta t$. As the mean correction is a function of the current time solution, we see that $\Delta t$ is not the primary factor in determining the best mean correction window. A larger timestep can mean that a local minimum at a smaller $\eta_C$ value is best, as is the case at $\epsilon = 0.1$ in Figure \ref{fig:RSWEs_meancor_err_over_etaC}.

\begin{figure}
     \centering
     \begin{subfigure}[b]{0.475\textwidth}
         \centering
         \includegraphics[width=\textwidth]{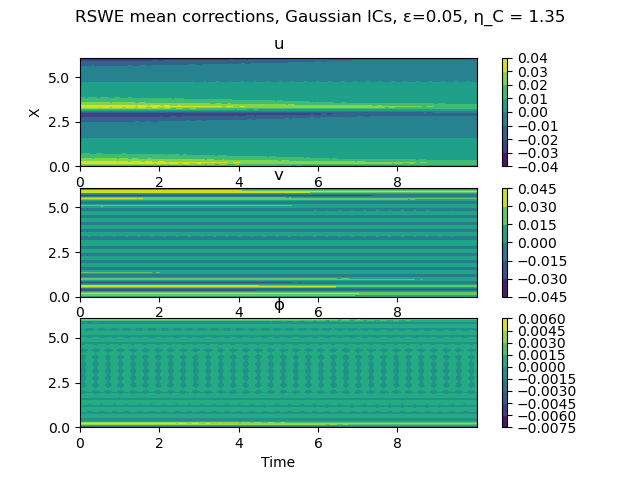}
         \caption{$\eta_C = 1.35$}
         \label{fig:meancor_visual_RSWE_Gaussian_etaC_1_35}
     \end{subfigure}
     \hfill
     \begin{subfigure}[b]{0.475\textwidth}
         \centering
    \includegraphics[width=\textwidth]{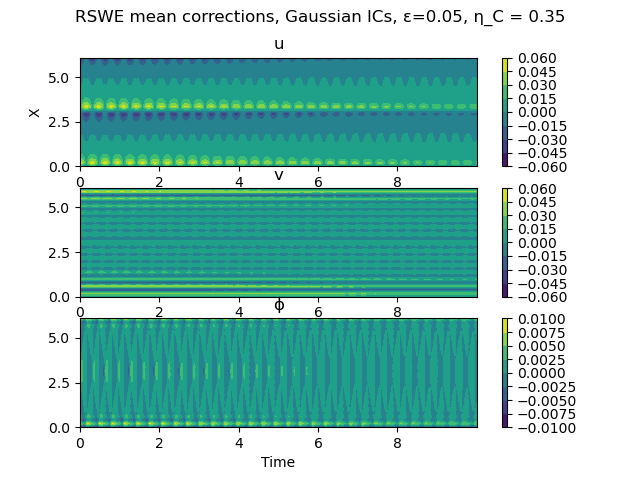}
         \caption{$\eta_C = 0.35$}
         \label{fig:meancor_visual_RSWE_Gaussian_etaC_0_35}
     \end{subfigure}
        \caption{Visualisations of each component of the RSWEs mean correction using two different $\eta_C$ values, at $\epsilon=0.05$. These are computed from $\left< \cdot \right>_r$ averages of a fine timestep $\vec{V}$ solution and transformed into physical space (analogously to Figure \ref{fig:KG_an_num_Cs}). The larger $\eta_C = 1.35$ in (a) is best for the smaller timestep of $\Delta t =0.1$; this corresponds to the third local minima in the error over $\eta_C$, like in Figure \ref{fig:RSWEs_meancor_err_over_etaC}. This larger window has a smoother evolution of the mean correction over time. The smaller $\eta_C = 0.35$ in (b) is best for the larger timestep of $\Delta t = 0.2$ and corresponds to the first minima in the error over $\eta_C$. This smaller averaging window allows for more local variation to be captured in $\vec{C}$. This is seen through shorter timescale oscillations in the $u$ component of (b) which are smoothed out in (a). The magnitude of each mean correction component reduces with a larger $\eta_C$, as there is more cancellation of oscillations.}
    \label{fig:meancor_visual_RSWE_Gaussian}
\end{figure}

The components of the mean correction in physical space are visualised in Figure \ref{fig:meancor_visual_RSWE_Gaussian} for two values of $\eta_C$ and $\epsilon = 0.05$. Similar mean dynamics are observed in both cases, with two periods across the spatial dimension in the $u$ component. This may be a result of the quadratic order of the nonlinearity halving the width of the initial Gaussian height perturbation; in the $u$ component we have $\mathcal{N}_u =  -u \frac{\partial u}{\partial x}$. The smaller window of $\eta_C = 0.35$ in Figure \ref{fig:meancor_visual_RSWE_Gaussian_etaC_1_35}, compared to $\eta_C = 1.35$ in Figure \ref{fig:meancor_visual_RSWE_Gaussian_etaC_0_35}, retains a larger range of frequencies in the local mean correction. This leads to more variation in $\vec{C}$ over time, especially in $C_u$, where there are visible oscillations with $\eta_C = 0.35$ but a much smoother time evolution with $\eta_C = 1.35$. There is a larger magnitude of each component of $\vec{C}$ with the smaller averaging window, as a larger $\eta_C$ leads to an increasing cancellation of the fast oscillations.  \par

\section{Discussion}
\label{sec:discussion}
This paper introduced a modified modulation variable mapping, with a local mean correction term, to improve phase-averaging accuracy in oscillatory, multiscale, differential equations. We provided algorithms for the local mean correction and a timestepping method that uses a finite averaging window; the inclusion of the local mean correction extends the standard finite phase-averaging framework of \cite{peddle_haut_wingate_2019}. Numerical examples with a classical mean correction in the swinging spring and KG-type equation show the viability of the mean correction beyond weakly nonlinear systems in \cite{Tao_simp_imp} into multiscale systems like equation (\ref{eq:standard_eq}). We saw accuracy improvements with a local mean correction in the one-dimensional RSWEs, which shows the potential for mean corrected phase-averaging in other weather and climate PDEs. \par
The local mean correction is a new quantity that combines the mean correction with the benefits of a finite phase-average. It retains low-frequency information in the mean correction, which is beneficial when using a finite phase-average of the timestepping equation. The best window length for the local mean correction, $\eta_C^*$, provides a balance between a sufficiently long interval to compute a mean value whilst being small enough to retain pertinent local variation. This is analogous to how the finite phase-average of the timestepping equation retains important low-frequency interactions and dynamics that are missing with a classical time-average \cite{peddle_haut_wingate_2019}.
\par
There is ongoing work on constructing a more efficient algorithm that precomputes a matrix exponential of $\exp(\pm \Delta t \mathcal{L})$ before the timestepping routine, as opposed to $\exp(\pm t \mathcal{L})$ at each timestep. Examples of this approach with phase-averaging include Lawson-RK methods, such as the fourth-order version used in \cite{hiroe_colin}, or a Strang-split method, as used in \cite{Haut_Wingate_2014}. Other possible additions to the locally mean corrected method include an optimisation algorithm to select $\eta^*$ and $\eta_C^*$ more efficiently, and extending the improved accuracy into slower oscillation regimes, i.e. $\epsilon \rightarrow 1$. Additionally, we would like to test the new mapping as a modified coarse solver in the oscillatory parareal method of \cite{peddle_haut_wingate_2019}. As the mean correction can be computed in parallel like the phase-averaging procedure, it could provide a cheap modification to the coarse solver to improve its accuracy. \par

\section*{Funding acknowledgements}
Timothy C. Andrews was supported by a grant from the Engineering and Physical Sciences Council (EPSRC) [grant number EP-V520317-1].  For the purpose of open access, the author has applied a ‘Creative Commons Attribution' (CC BY) licence to any author accepted manuscript version arising from this submission.

\section*{Appendix A: Cancellation of oscillations example for the swinging spring $\vec{C}_{\infty}$}
\label{appendix:canc_of_osc}

Consider the first term in the $x$ component of $\mathcal{N}(e^{-t \mathcal{L} }\vec{W})$ in \ref{eq:swing_spring_N_w},
\begin{equation}
    \frac{i \lambda}{4 \omega_R} A \tilde{x} \tilde{z} =  \frac{i \lambda}{4 \omega_R} e^{-i (\rho + 1) \omega_R t} \tilde{x}(t) \tilde{z}(t).
\end{equation}

\noindent Applying (\ref{eq:C_inf}) to obtain the mean value over all time,
\begin{equation}
    \lim_{T_C \rightarrow \infty } \frac{1}{T_C} \int_{r=0}^{T_C} \frac{i \lambda}{4 \omega_R} e^{-i (\rho + 1) \omega_R (t + r)} \tilde{x}(t) \tilde{z}(t) dr.
\end{equation} 

\noindent Phase-averaging is performed at a fixed solution value, allowing most terms to be taken outside of the integral,
\begin{equation}
    \frac{i \lambda}{4 \omega_R} e^{-i (\rho + 1) \omega_R t} \tilde{x}(t) \tilde{z}(t)  \lim_{T_C \rightarrow \infty } \frac{1}{T_C} \int_0^{T_C}  e^{-i (\rho + 1) \omega_R r} dr.
\end{equation}

\noindent This integral evaluates to zero by the following property of complex exponentials,
\begin{equation}
    \lim_{T_C \rightarrow \infty } \frac{1}{T_C} \int_0^{T_C} e^{i \phi t} dt = 
\begin{cases}
    0, if \ \phi \neq 0, \\
    1, if \ \phi = 0. 
\end{cases}
\end{equation}
\par
Every term in (\ref{eq:swing_spring_N_w}) with a complex exponential component ($A, B, C, D, E, F$) will have this cancellation of oscillations, resulting in a zero mean value and no contribution to $\vec{C}_{\infty}$.  

\section*{Appendix B: KG-type equation classical mean correction}
\label{appendix:kg_C_an}

We begin by expressing the KG-type PDE nonlinear operator (\ref{eq:KG_operators}), $\mathcal{N} = [\mathcal{N}_a, \mathcal{N}_b]^{\text{T}}$, in terms of the mean corrected modulation variable in Fourier space, $\vec{\hat{W}} = \left[ \hat{c},\hat{d} \right]^{\text{T}}$, using a mapping of the matrix exponential of (\ref{eq:KG_matexp}),
\begin{equation}
    \mathcal{N} \left( e^{ - \frac{t \mathcal{L}}{\epsilon}} \vec{\hat{W}} \right) 
=
\begin{bmatrix}
0\\
- \left( \hat{c} \frac{\epsilon}{\omega} \cos \left(\frac{\omega t}{\epsilon}\right) + \hat{d} \frac{\epsilon}{\omega}  \sin \left(\frac{\omega t}{\epsilon}\right) \right) \circledast \left( \hat{c} \frac{\epsilon}{\omega} \cos \left(\frac{\omega t}{\epsilon}\right) + \hat{d} \frac{\epsilon}{\omega}  \sin \left(\frac{\omega t}{\epsilon}\right) \right)  \\
\end{bmatrix},
\label{eq:KG_N_w}
\end{equation}

\noindent where $\circledast$ is a spatial circular convolution, defined in its discrete form for any two length $N$ arrays, $a$ and $b$, as
\begin{equation}
    g = a \circledast b \ \ \ \ \rightarrow g_i = \frac{1}{N} \sum_{j=0}^{N-1} a_j b_{(i - j) \ \mathrm{mod} \ N},
    \label{eq:conv_sum}
\end{equation}

\noindent with subscripts denoting the entry of a discrete array. \par
We can expand the non-zero nonlinear component as $\mathcal{N}_b = \mathcal{N}^I + \mathcal{N}^{II} + \mathcal{N}^{III} + \mathcal{N}^{IV}$, as $\circledast$ distributes over addition,
\begin{subequations}
\begin{align}
    \mathcal{N}^I &= - \left( \hat{c} \frac{\epsilon}{\omega} \cos \left(\frac{\omega t}{\epsilon}\right) \right) \circledast \left( \hat{c} \frac{\epsilon}{\omega} \cos \left(\frac{\omega t}{\epsilon}\right) \right),
     \\
    \mathcal{N}^{II} &= - \left( \hat{c} \frac{\epsilon}{\omega} \cos \left(\frac{\omega t}{\epsilon}\right) \right) \circledast \left( \hat{d} \frac{\epsilon}{\omega} \sin \left(\frac{\omega t}{\epsilon}\right) \right) ,
     \\
    \mathcal{N}^{III} &= - \left( \hat{d} \frac{\epsilon}{\omega} \sin \left(\frac{\omega t}{\epsilon}\right) \right) \circledast \left( \hat{c} \frac{\epsilon}{\omega} \cos \left(\frac{\omega t}{\epsilon}\right) \right),
     \\
    \mathcal{N}^{IV} &= - \left( \hat{d} \frac{\epsilon}{\omega} \sin \left(\frac{\omega t}{\epsilon}\right) \right) \circledast \left( \hat{d} \frac{\epsilon}{\omega} \sin \left(\frac{\omega t}{\epsilon}\right) \right). 
\end{align}
\label{eq:KG_N_convs}
\end{subequations}

The mean correction can be decomposed into averages of each component of $\mathcal{N}_b$, using the notation of $C^i = \left< \mathcal{N}^i \right>_r$. For a classical mean correction, we evaluate each term in (\ref{eq:KG_N_convs}) in the limit of $T_C \rightarrow \infty$ by applying (\ref{eq:C_inf}). We then use the following results of integrating a product of trigonometric functions,

\begin{subequations}
\begin{align}
\lim_{T_C \to \infty} \frac{1}{T_C}  \int_0^{T_C}
 \cos \left(\frac{\omega_a t}{\epsilon}\right) \cos \left(\frac{\omega_b t}{\epsilon}\right) 
&=
\begin{cases}
\frac{1}{2}, ~\text{if } \omega_a = \omega_b, \\
0, \text{if } \omega_a \ \neq \omega_b,
\end{cases}
\label{eq:trig_int_cos_cos}
\\
        \lim_{T_C \to \infty} \frac{1}{T_C} \int_0^{T_C}
\cos \left(\frac{\omega_a t}{\epsilon}\right) \sin \left(\frac{\omega_b t}{\epsilon}\right)
&=
0, \forall ~\omega_a, \omega_b ,
\label{eq:trig_int_cos_sin}
\\
    \lim_{T_C \to \infty} \frac{1}{T_C} \int_0^{T_C}
\sin \left(\frac{\omega_a t}{\epsilon}\right) \cos \left(\frac{\omega_b t}{\epsilon}\right)
&=
0, \forall ~\omega_a, \omega_b,
\label{eq:trig_int_sin_cos}
\\
\lim_{T_C \to \infty} \frac{1}{T_C}  \int_0^{T_C}
 \sin \left(\frac{\omega_a t}{\epsilon}\right) \sin \left(\frac{\omega_b t}{\epsilon}\right) 
&=
\begin{cases}
\frac{1}{2}, ~\text{if } \omega_a = \omega_b, \\
- \frac{1}{2}, ~\text{if } \omega_a = - \omega_b, \\
0, \text{if } |\omega_a| \ \neq |\omega_b|.
\end{cases}
\label{eq:trig_int_sin_sin}
\end{align}
\end{subequations}

By (\ref{eq:trig_int_cos_sin}) and (\ref{eq:trig_int_sin_cos}), $C^{II}$ and $ C^{III}$ will have a mean value of zero in the $T_C \rightarrow \infty$ limit. $C^{I}$ and $ C^{IV}$ will have a non-zero contribution in the classical mean correction where $|\omega_a| = |\omega_b|$. The KG-type dispersion relation of $\omega(k) = \sqrt{1 + k^2}$ is a positive definite and even function, as $\omega(k) = \omega(-k)$. This means that two wavenumbers of equal magnitude will have the same frequencies and we cannot have $\omega_a = -\omega_b$. \par
Using (\ref{eq:trig_int_cos_cos}) and (\ref{eq:trig_int_sin_sin}), we are left with a classical mean correction of $\vec{C}_{\infty} = \left[ 0, C^b \right]^{\text{T}}, C^b = C^I + C^{IV}$, where

\begin{equation}
    C^I_i = 
\begin{dcases}
\frac{1}{2 N_x} \sum_{j=0}^{N_x-1} \frac{ \hat{c}_j \hat{c}_{-j \ \mathrm{mod} \ N} }{ \omega_j^2} ,  &\text{for } i=0, \\
\frac{1}{2 N_x} \left( \frac{ \hat{c}_{i/2}^2 }{ \omega_{i/2}^2} + \frac{ \hat{c}_{(i+N_x)/2}^2 }{ \omega_{(i+N_x)/2}^2} \right), &\text{for } i ~\text{even}, \\
0, &\text{for} \ i ~\text{odd},
\end{dcases}
\end{equation}

\begin{equation}
C^{IV}_i = 
\begin{dcases}
\frac{1}{2 N_x} \sum_{j=0}^{N_x-1} \frac{ \hat{d}_j \hat{d}_{-j \ \mathrm{mod} \ N} }{ \omega_j^2} ,  &\text{for} \ i=0, \\
\frac{1}{2 N_x} \left( \frac{ \hat{d}_{i/2}^2 }{ \omega_{i/2}^2} + \frac{ \hat{d}_{(i+N_x)/2}^2 }{ \omega_{(i+N_x)/2}^2} \right), &\text{for} \ i ~\text{even}, \\
0, &\text{for} \ i ~\text{odd}.
\end{dcases}
\end{equation}

\printbibliography

\end{document}